\documentclass[12pt]{amsart}
\usepackage{amsmath}
\usepackage{amssymb}
\usepackage{latexsym}
\newcommand{\halmos}{\rule{5pt}{5pt}}

\newcommand{\Heun}{Heun~equation}

\numberwithin{equation}{section}
\newtheorem{thm}{Theorem}[section]
\newtheorem{prop}[thm]{Proposition}

\theoremstyle{definition} 

\theoremstyle{remark}

\begin{document}
\title[Degenerations of $q$-Heun equation]
{Degenerations of $q$-Heun equation}
\author{Chihiro Sato}
\address{Department of Mathematics, Ochanomizu University, 2-1-1 Otsuka, Bunkyo-ku, Tokyo 112-8610, Japan}
\email{clover.chihiro@icloud.com}
\author{Kouichi Takemura}
\address{Department of Mathematics, Ochanomizu University, 2-1-1 Otsuka, Bunkyo-ku, Tokyo 112-8610, Japan}
\email{takemura.kouichi@ocha.ac.jp}
\subjclass[2010]{33E10,39A13}
\keywords{Heun's differential equation, $q$-Heun equation, degeneration}
\begin{abstract}
We obtain several degenerations of the $q$-Heun equation by considering the linear $q$-difference equations associated to several $q$-Painlev\'e equations.
We establish definitions of the confluent $q$-Heun equation, the biconfluent $q$-Heun equation and the doubly confluent $q$-Heun equation, and investigate limit procedures to the corresponding differential equations. 
%For more discovery of series solutions of the Heun equation and its degeneration, we investigate degeneration of the q-Heun equation which is a q-difference equation of the Heun equation.
\end{abstract}
\maketitle

\section{Introduction}

The $q$-Heun equation is a $q$-difference analogue of the Heun equation, which is given as the form
     \begin{equation}
        (\alpha_2x^2+\alpha_1x+\alpha_0)f(qx)
        -(\beta_2x^2+\beta_1x+\beta_0)f(x)
        +(\gamma_2x^2+\gamma_1x+\gamma_0)f(x/q)=0,        
      \end{equation}
under the condition $\alpha_2\alpha_0\gamma_2\gamma_0\not=0$.
It was introduced by Hahn \cite{Hahn} in 1971, and was rediscovered by the second author \cite{TakR} in a study of degenerations of Ruijsenaars-van Diejen operator \cite{vD0,RuiN}.
Recall that the Heun equation or Heun's differential equation is a standard form of the second order linear differential equation with four regular singularities on the Riemann sphere, and it is written as 
\begin{equation}
  \frac{d^2y}{dz^2}
  +\left(\frac{\gamma}{z}+\frac{\delta}{z-1}+\frac{\epsilon}{z-t}\right)\frac{dy}{dz}
  +\frac{\alpha\beta z-B}{z(z-1)(z-t)}y=0,
\label{eq:Heun}
\end{equation}
under the condition $ \gamma+\delta+\epsilon=\alpha+\beta+1$.
The parameter $B$ is an accessory parameter, which is independent from the local exponents.
It is known that Heun's differential equation is recovered from the $q$-Heun equation as $q \to 1$ (\cite{Hahn,TakR}).
On the $q$-Heun equation, we may regard $\beta _1$ as the accessory parameter.

The Heun equation (HE) admits degenerations to confluent Heun equation (CHE), biconfluent Heun equation (BHE), doubly confluent Heun equation (DHE) and triconflunent Heun equation (THE).
\begin{align}
& \mbox{CHE: } \; \frac{d^2y}{dz^2} +\left(\frac{\gamma}{z}+\frac{\delta}{z-1}-\beta\right)\frac{dy}{dz}
      +\frac{-\alpha \beta z+B}{z(z-1)}y=0 , \label{eq:CHE} \\
& \mbox{BHE: } \; \frac{d^2y}{dz^2} +\left(-z-\delta+\frac{\gamma}{z}\right)\frac{dy}{dz} +\left(-\alpha+\frac{B}{z}\right)y=0 , \label{eq:BHE}\\
& \mbox{DHE: } \; \frac{d^2y}{dz^2} +\left(-1-\frac{\gamma}{z}-\frac{\delta}{z^2}\right)\frac{dy}{dz} +\left(-\frac{\alpha}{z}+\frac{B}{z^2}\right)y=0 , \label{eq:DHE} \\
& \mbox{THE: } \; \frac{d^2y}{dz^2} +\left(-z^2-\gamma\right)\frac{dy}{dz} +\left(\alpha z+B\right)y=0 . \label{eq:THE}
\end{align}
The degeneration scheme is given as follows.
\begin{center}
\begin{tabular}{ccccccccc}
    \Heun & $\! \rightarrow \! $ & Confluent HE & $\rightarrow \!$ & Biconfluent HE & &&\\
    (HE)  &               & (CHE)        &  $\searrow \! $   & (BHE)          &$ \! \searrow \! $&&\\
          &&&& Doubly confluent HE & $\! \rightarrow \! $ & Triconfluent HE &\\
          &&&& (DHE)               &               & (THE)           &\\
\end{tabular}
\end{center}\vspace{0cm}
On the Heun equation (\ref{eq:Heun}), all the singularities $\{ 0,1,t,\infty \}$ are regular.
The confluent Heun equation (CHE) has two regular singularities $\{ 0,1 \}$ and an irregular singularity $\{ \infty \}$.
We may regard that the irregular singularity $z = \infty $ is obtained by confluencing the singularity $z=t$ of the Heun equation to $z=\infty $.
The Heun equation and the confluent ones sometimes appear in physics, especially in the analysis of black hole associated with the Kerr solution of the general relatibity \cite{Htd,NM}.

It is known that the Heun equation is related to the Painlev\'e equation (see \cite{TakHP,TakM}).
The $q$-Heun equation is also related to the $q$-Painleve equation.
Jimbo and Sakai \cite{JS} introduced a $q$-analogue of the sixth Painlev\'e equation as
\begin{equation}
f\overline{f} = \nu_3\nu_4\, \frac{ ( g - \nu_5 /\kappa_2 ) ( g - \nu_6/\kappa_2 )}{ ( g - 1/\nu_1 ) ( g - 1/\nu_2) }, \quad
g\underline{g} = \frac{1}{\nu_1\nu_2}\frac{( f - \kappa_1 /\nu_7 ) ( f - \kappa_1/\nu_8 )}{(f - \nu_3)(f - \nu_4)}. \label{eq:qPD15}
\end{equation}
Here, $\overline{f} $ (resp. $\underline{f}$) is the consequence of the time evolution $t \mapsto qt$ (resp. the inverse time evolution $t \mapsto t/q$) to $f$.
The time evolution for the parameters is given by $\overline{\kappa }_1= \kappa_{1} /q$, $\overline{\kappa }_2= q \kappa_{2}$ and $ \overline{\nu }_i= \nu _{i}$ $(i=1,2,\dots ,8)$, and the parameters are constrained by the relation $\kappa_{1} ^2 \kappa_{2} ^2 =q \nu _1 \nu _2 \dots \nu _8$. 
%Note that the expression of equation (\ref{eq:qPD15}) is due to Kajiwara, Noumi and Yamada written in the surver \cite{KNY}.
Jimbo and Sakai obtained equation (\ref{eq:qPD15}) by considering the connection preserving deformation of the $q$-linear difference equations written as $Y(qx,t)=A(x,t)Y(x,t)$, where $A(x,t) $ is a certain $2\times 2$ matrix, and the connection preserving deformation is equivalent to a compatibility condition for a Lax form.
Sakai provided a list of the second order discrete Painlev\'e equations in \cite{Sak}.
Each member of the $q$-difference Painlev\'e equation was essentially labelled by the affine root system from the space of initial conditions.
%Sakai \cite{Sak} realized the time evolution of the discrete Painlev\'e equations by the symmetry of the Weyl group, which originates from the Cremona action on a family of surfaces which are related with the space of initial conditions.
%Equation (\ref{eq:qPD15}) has the symmetry of the extended affine Weyl group $D^{(1)}_5$, and the equation is denoted by $q$-$P(D^{(1)}_5)$.
Murata \cite{Mrt} investigated degenerations of the Lax forms from the ones of Jimbo and Sakai.
He used the following table to express the degenerations of the $q$-Painlev\'e equations.
\\
\begin{center}
  \begin{tabular}{ccccccccccccccc}
%        [M]&
        $A_0^*$&
        $\! \rightarrow \! $&$A_1$&
        $\! \rightarrow \! $&$A_2$&
        $\! \rightarrow \! $&$A_3$&
        $\! \rightarrow \! $&$A_4$&
        $\! \rightarrow \! $&$A_5$&
        $\! \rightarrow \! $&$A_6$&
        $\! \rightarrow \! $&$A_7$\\
        &&&&&&&&&$\! \searrow \! $&&$ \! \nearrow \!$&&$ \! \nearrow \! $\\
        &&&&&&&&&&$A_5^\#$&$\rightarrow$&$A_6^\#$&$\rightarrow$&$A_7'$
  \end{tabular}
\end{center}
The sixth $q$-Painlev\'e equation of Jimbo and Sakai corresponds to type $A_3$, and Murata obtained the Lax forms for the $q$-Painlev\'e equation of type $A_4$ by considering the degeneration from type $A_3$.
He also obtained the Lax forms for the $q$-Painlev\'e equation for the decendants from type $A_4$.
The Lax forms were described as the pair of the systems of linear $q$-difference equations $Y(qx,t)=A(x,t)Y(x,t)$ and $Y(x,qt)=B(x,t)Y(x,t)$ where $Y(x,t)= \! ^t (y_1(x,t), y_2(x,t) ) $ and $A(x,t)$ and $B(x,t)$ are $2\times 2$ matrices whose coefficients contain the variables $\lambda $ and $\mu $, and the $q$-Painlev\'e equations with the variables $\lambda $ and $\mu $ were described as the compativility condition $A(x,qt)B(x,t)=B(qx,t)A(x,t)$ for these two systems.
Yamada obtained the Lax pair of the elliptic Painlev\'e equation in \cite{Ye}, which is located at the top on degerenations, and he also obtained the Lax pairs for the $q$-Painlev\'e equations in \cite{Y}.
%Yamada studied the Lax pairs of the discrete Painlev\'e equations extensively \cite{Ye,Y}, and a list of the Lax pairs was provided in the survey \cite{KNY} by Kajiwara, Noumi and Yamada.
Kajiwara, Noumi and Yamada \cite{KNY} wrote a survey of discrete Painlev\'e equations, and they gave a list of the discrete Painlev\'e equations including the $q$-Painlev\'e equations, the Lax pairs associated to each discrete Painlev\'e equation and the extended affine Weyl group symmetries and so on.
They used the labelling of the $q$-Painlev\'e equations by the symmetry as the following table.
\\
\begin{center}
{
  \begin{tabular}{ccccccccccccccc}
%       {[KNY]}&
        $E_8^{(1)}$&
        $\! \rightarrow \! $&$E_7^{(1)}$&
        $\! \rightarrow \! $&$E_6^{(1)}$&
        $\! \rightarrow \! $&$D_5^{(1)}$&
        $\! \rightarrow \! $&$A_4^{(1)}$&
        $\! \rightarrow \! $&$E_3^{(1)}(b)$&
        $\! \rightarrow \! $&$E_2^{(1)}(b)$&
        $\! \rightarrow \! $&$A_1^{(1)}$\\
        &&&&&&&&&$\! \searrow \! $&&$ \! \searrow \! $&&$\! \searrow \! $\\
        &&&&&&&&&&$E_3^{(1)}(a)$&
        $\rightarrow$&$E_2^{(1)}(a)$&
        $\rightarrow$&$\underset{|\alpha^2|=8}{A_1^{(1)}}$
  \end{tabular}}
\end{center}
The sixth $q$-Painlev\'e equation of Jimbo and Sakai corresponds to type $D^{(1)}_5$, and the equation is denoted by $q$-$P(D^{(1)}_5)$ in \cite{KNY}.
A Lax pair of $q$-$P(D^{(1)}_5)$ given in \cite{KNY} is written as
\begin{align}
L_{1} &= \Bigl\{\frac{z(g\nu_{1}-1)(g\nu_{2}-1)}{qg}-\frac{\nu_{1}\nu_{2}\nu_{3}\nu_{4}(g- \nu_{5}/\kappa_{2})(g- \nu_{6}/\kappa_{2})}{fg}\Bigr\} \label{eq:D5L1} \\
&+\frac{\nu_{1}\nu_{2}(z- q \nu_{3})(z - q \nu_{4})}{q(qf-z)}(g-T^{-1}_{z})+\frac{(z- \kappa_{1} /\nu_{7})(z - \kappa_{1}/\nu_{8})}{q(f-z)}\Bigl(\frac{1}{g} - T_{z}\Bigr) , \nonumber \\
L_{2} &= \Bigl(1-\frac{f}{z} \Bigr)T+T_{z}-\frac{1}{g} . \nonumber 
\end{align}
Here $T_{z}$ represents the transformation $z\mapsto qz$ and $T$ represents the time evolution such as $T(f)=\overline{f}$.
It follows from the compatibility condition for the Lax operators $L_1$ and $L_2$ (see \cite{KNY}) that the parameters $f$ and $g$ satisfies the $q$-Painlev\'e equation $q$-$P(D^{(1)}_5)$ given in equation (\ref{eq:qPD15}).

We now explain a relationship between the equation $q$-$P(D^{(1)}_5)$ and the $q$-Heun equation (\cite{STT2,TakR}).
We substitute $f= \nu_{4} $ to the equation $L_1 y(z)=0 $.
Then, we have
\begin{align}
    &(z-q \nu_{3})(z- \nu_{4})y(z/q)+\frac{1}{\nu_{1}\nu_{2}}\Bigl(z-\frac{\kappa_{1}}{\nu_{7}}\Bigr)\Bigl(z-\frac{\kappa_{1}}{\nu_{8}}\Bigr)y(qz)  \\
    &-\Bigl[\Bigl(\frac{1}{\nu_{1}}+\frac{1}{\nu_{2}}\Bigr)z^{2} - \Bigl( \frac{ (\nu _3 \nu_7 - \kappa_{1} )(\nu _4 \nu_7 - \kappa_{1} )}{\nu_{1} \nu_{2} \nu_{4} \nu_{7} \nu_{8} }\frac{1}{g} +\frac{\nu_{4}}{\nu_{1}} + \frac{\nu_{4}}{\nu_{2}}+\frac{q\nu_{3}\nu_{5}}{\kappa_{2}}+\frac{q\nu_{3}\nu_{6}}{\kappa_{2}}\Bigr)z \nonumber \\
& \qquad +\frac{q \nu_{3}\nu_{4}{\nu_{5}}^{1/2}{\nu_{6}}^{1/2}}{\kappa_{2}}\Bigl\{\Bigl(\frac{\nu_{6}}{\nu_{5}}\Bigr)^{1/2}+\Bigl(\frac{\nu_{6}}{\nu_{5}}\Bigr)^{-1/2}\Bigr\}
    \Bigr] y(z) =0 . \nonumber
\end{align}
Hence, we obtain the $q$-Heun equation, and we may regard $g$ as an accessory parameter.

In this paper, we obtain degenerate $q$-Heun equations from the linear equations $L_1y(z)=0$ which are associated to the $q$-Painlev\'e equation $q$-P($A_4^{(1)}$) and its descendants.
%As a result, we obtain several degenerations of the $q$-Heun equation by considering the linear $q$-difference equations associated to several $q$-Painlev\'e equations.
Based on partial members of the degenerate $q$-Heun equations, we establish definitions of the confluent $q$-Heun equation, the biconfluent $q$-Heun equation and the doubly confluent $q$-Heun equation, and investigate limit procedures to the corresponding differential equations. 

This paper is organized as follows.
In Section \ref{sec:gt}, we review results on gauge transformation, which are used in this paper.
We recall the Newton polygon for the $q$-difference equation, which visualizes features of the $q$-difference equation.
In Section \ref{sec:Mrt}, we obtain degenerations of $q$-Heun equations from the Lax forms of Murata.
In Section \ref{sec:KNY}, we obtain degenerations of $q$-Heun equations from the Lax pairs of Kajiwara, Noumi and Yamada.
In Section \ref{sec:cqHs}, we introduce definitions of the confluent $q$-Heun equation, the biconfluent $q$-Heun equation and the doubly confluent $q$-Heun equation, and investigate limit procedures to the corresponding differential equations.
In Section \ref{sec:rmk}, we give concluding remarks.

\section{Newton polygon and gauge transformation} \label{sec:gt}

We investigate the $q$-difference equation which is written in the form
     \begin{gather}
        (\alpha_2x^2+\alpha_1x+\alpha_0)f(qx)
        -(\beta_2x^2+\beta_1x+\beta_0)f(x)
        +(\gamma_2x^2+\gamma_1x+\gamma_0)f(x/q)=0.       
\label{eq:qHeun0}     
     \end{gather}
We associate the $q$-difference equation (\ref{eq:qHeun0}) with the Newton diagram as \cite{Sau} (see also \cite{Adams31,Ohy}).
In the Newton diagram, the black circle means a coefficient which is not equal to zero, and the write circle means a coefficient which is equal to zero.
We show examples of the Newton diagrams.
The non-degenerate $q$-Heun equation is the one which satisfies the condition $\alpha_2\alpha_0\gamma_2\gamma_0\not=0$, and the corresponding Newton diagram is written as follows.
\begin{center}
\begin{picture}(200,85)(0,0)
%  \put(15,172){q-P($A_4$)}
    \put(70,60){$\alpha _2$}
    \put(70,40){$\alpha _1$}
    \put(70,20){$\alpha _0$}
    \put(50,60){$\beta _2$}
    \put(50,40){$\beta _1$}
    \put(50,20){$\beta _0$}
    \put(30,60){$\gamma _2$}
    \put(30,40){$\gamma _1$}
    \put(30,20){$\gamma _0$}
% \put(110,72){q-P($A_5^\#$)}
    %\put(90,110){\vector(1,-4){15}}
    \put(125,60){\circle*{5}}
    \put(142,56){$*$}
    \put(165,60){\circle*{5}}
%    \put(125,40){\circle*{5}}
%    \put(145,40){\circle*{5}}
%    \put(165,40){\circle*{5}}
    \put(122,36){$*$}
    \put(142,36){$*$}
    \put(162,36){$*$}
    \put(125,20){\circle*{5}}
    \put(142,16){$*$}
%    \put(145,20){\circle*{5}}
    \put(165,20){\circle*{5}}
    \qbezier(125,60)(135,60)(145,60)
    \qbezier(145,60)(155,60)(165,60)
%    \qbezier(165,40)(165,30)(165,20)
    \qbezier(125,20)(135,20)(145,20)
    \qbezier(165,20)(155,20)(145,20)
%    \qbezier(145,20)(135,30)(125,40)
  \end{picture}
\end{center}
The asterisk $(*)$ means a coefficient which can be either zero or non-zero.
On the other hand, we call the $q$-difference equation (\ref{eq:qHeun0}) which satisfies the condition $\alpha_2\alpha_0\gamma_2\gamma_0 =0$ by the degenerate $q$-Heun equation.
The $q$-difference equation (\ref{eq:qHeun0}) which satisfies the condition $\alpha_2= \beta _2 = \gamma_2 =0$ is called the hypergeometric type (see \cite{Ohy}).

We remind some formulas for gauge transformation of $q$-difference equation.
\begin{prop} \label{prop:GT} $ $\\
(i) If $y(x)$ satisfies $q^{-\lambda}a(x)g(x/q)+b(x)g(x)+q^{\lambda}c(x)g(qx)=0$, then the function $h(x)=x^{\lambda}y(x)$ satisfies $a(x)g(x/q)+b(x)g(x)+c(x)g(qx)=0$.\\
(ii) If $y(x)$ satisfies $(1-\alpha x)a(x)g(x/q)+b(x)g(x)+c(x)g(qx)=0$, then the function $u(x)=(q\alpha x;q)_{\infty}y(x)$ satisfies $a(x)g(x/q)+b(x)g(x)+(1-q \alpha x)c(x)g(qx)=0$.
(iii) Set $\theta _q (x)= (q,-x,-q/x;q)_{\infty }$.
If $y(x)$ satisfies $\alpha x a(x)g(x/q)+b(x)g(x)+c(x)g(qx)=0$, then the function $w(x)=\theta _q (q \alpha x) y(x)$ satisfies $a(x)g(x/q)+b(x)g(x)+q \alpha x c(x)g(qx)=0$.
\end{prop}
See \cite{Ohy,STT1,MST}.
Note that  $\theta _q (x)$ satisfies $x \theta _q (xq)= \theta _q (x) $.

On the equation
\begin{align}
& f(qx) -(\beta_1x+\beta_0)f(x) +\gamma (1- \tilde{\gamma }_1 x)( 1 - \tilde{\gamma }_2 x)f(x/q)=0, \quad \gamma \tilde{\gamma }_1 \tilde{\gamma }_2 \neq 0 , 
\label{eq:a0g20}
\end{align}
we apply the gauge transformation $f(x) = (q\tilde{\gamma }_1 x;q)_{\infty} g(x)$.
Then, we obtain
\begin{align}
&  (1- q \tilde{\gamma }_1 x) g(qx) -(\beta_1x+\beta_0)g (x) +\gamma( 1 - \tilde{\gamma }_2 x)f(x/q)=0, \quad \gamma \tilde{\gamma }_1 \tilde{\gamma }_2 \neq 0 , 
\label{eq:a10c10}
\end{align}
by Proposition \ref{prop:GT} (ii), and it is the hypergeometric type.
The Newton diagrams of these two equations are written as follows.
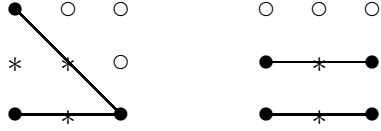
\begin{center}
\begin{picture}(400,85)(0,0)
    \put(165,60){\circle{5}}
    \put(145,60){\circle{5}}
    \put(125,60){\circle*{5}}
%    \put(165,40){\circle*{5}}
%    \put(145,40){\circle*{5}}
%    \put(125,40){\circle*{5}}
    \put(165,40){\circle{5}}
    \put(142,36){$*$}
    \put(122,36){$*$}
    \put(165,20){\circle*{5}}
    \put(142,16){$*$}
%    \put(145,20){\circle*{5}}
    \put(125,20){\circle*{5}}
    \qbezier(165,20)(155,30)(145,40)
    \qbezier(145,40)(135,50)(125,60)
%    \qbezier(165,40)(165,30)(165,20)
    \qbezier(125,20)(135,20)(145,20)
    \qbezier(165,20)(155,20)(145,20)
    \put(220,60){\circle{5}}
    \put(240,60){\circle{5}}
    \put(260,60){\circle{5}}
    \put(220,40){\circle*{5}}
    \put(237,36){$*$}
%    \put(240,40){\circle*{5}}
    \put(260,40){\circle*{5}}
    \put(220,20){\circle*{5}}
    \put(237,16){$*$}
%    \put(240,20){\circle{5}}
    \put(260,20){\circle*{5}}
    \qbezier(220,40)(240,40)(260,40)
    \qbezier(220,20)(240,20)(260,20)
  \end{picture}
\end{center}
We do not regard equation (\ref{eq:a0g20}) as the $q$-Heun type.

\section{Degenerations of $q$-Heun equations arising from the list of Murata} \label{sec:Mrt}

Murata (\cite{Mrt}) investigated degenerations of the Lax forms of several $q$-Painlev\'e equations.
The Lax forms were described as the pair of the systems of linear $q$-difference equations $Y(qx,t)=A(x,t)Y(x,t)$ and $Y(x,qt)=B(x,t)Y(x,t)$ where $Y(x,t):={^t\left(y_1(x,t), y_2(x,t)\right)}$ and $A(x,t)$ and $B(x,t)$ are $2\times 2$ matrices whose coefficients contain the variables $\lambda $ and $\mu $, and the $q$-Painlev\'e equations with the variables $\lambda $ and $\mu $ were described as the compativility condition $A(x,qt)B(x,t)=B(qx,t)A(x,t)$ for these two systems.

In this section, we obtain degenerations of $q$-Heun equations from the linear $q$-difference equation $Y(qx,t)=A(x,t)Y(x,t)$ by restricting the values $(\lambda , \mu ) $.
Note that non-degenerate $q$-Heun equation was obtained in \cite{TakR} from the linear $q$-difference equation $Y(qx,t)=A(x,t)Y(x,t)$ associated to the $q$-Painlev\'e equtaion $q$-P($A_3$).
The equation $q$-P($A_4$) and its decendants are obtained by degeneration from $q$-P($A_3$).
Therefore, we investigate the case $q$-P($A_4$) and its descendants.

\subsection{$q$-P($A_4$)}
Let $\theta_1, \theta_2, \kappa_1, \kappa_2, a_1, a_2, a_3$ be the parameters which satisfy $\theta_1\theta_2=-\kappa_1\kappa_2a_1a_2a_3$.
On the equation $Y(qx,t)=A(x,t)Y(x,t)$ for $q$-P($A_4$) by Murata, the matrix $A(x,t) $ is defined as 
\begin{align}
& A(x,t)=\begin{pmatrix}\kappa_1\{(x-\lambda)(x-\alpha)+\mu_1\}&\omega(x-\lambda)\\
          \kappa_1\omega^{-1}(\gamma x+\delta)&\kappa_2(x-\lambda)+\mu_2\end{pmatrix},\\
&  \mu_1=\frac{(\lambda-a_1t)(\lambda-a_2t)}{q\kappa_1\mu}, \;
  \mu_2=q\kappa_1\kappa_2\mu(\lambda-a_3), \nonumber \\
&  \alpha=\frac{\lambda^{-1}\left\{(\theta_1+\theta_2)t-\kappa_1\mu_1-\mu_2\right\}+\kappa_2}
              {\kappa_1}, \;
  \gamma=\mu_2-\kappa_2\left\{(2\lambda+\alpha)-(a_1+a_2)t-a_3\right\}, \nonumber \\
&  \delta=-\lambda^{-1}
          \left\{\kappa_2a_1a_2a_3t^2-(\alpha\lambda+\mu_1)(\kappa_2\lambda-\mu_2)\right\}. \nonumber
\end{align}
Note that it follows that $\det A(x,t)=\kappa_1\kappa_2(x-a_1t)(x-a_2t)(x-a_3)$ and the eigenvalues of the matrix $A(0,t)$ are $\theta_1t,\theta_2t$.
The parameters $\lambda $ and $\mu $ are the dependent variables of the $q$-Painlev\'e equation with respect to the independent varible $t$, though they are fixed in the current setting.
Set $Y(x,t)= \! ^t (y_1(x), y_2(x))$.
It follows from the equation $Y(qx,t)=A(x,t)Y(x,t)$ that the function $y_1(x) $ satisfies the single $q$-difference equation
\begin{align}
& y_1(q^2x)
  -\frac{q^2\kappa_1x^3-b_2x^2+b_1x-\lambda(\theta_1+\theta_2)t}{x-\lambda}y_1(qx) \\
& \qquad \qquad  +\frac{\kappa_1\kappa_2(qx-\lambda)(x-a_1t)(x-a_2t)(x-a_3)}{x-\lambda}y_1(x)=0, \nonumber
\end{align}
where
\begin{align}
&  b_2=q(q+1)\kappa_1\lambda
      -\frac{q^2\kappa_1\kappa_2\mu(\lambda-a_3)}{\lambda}
      +2q\kappa_2+\frac{q(\theta_1+\theta_2)t}{\lambda}
      -\frac{(\lambda-a_1t)(\lambda-a_2t)}{\mu\lambda},\\
&  b_1=q\kappa_1\lambda^2-q\kappa_1\kappa_2\mu(\lambda-a_3)
      -2q\kappa_2\lambda+(q+1)(\theta_1+\theta_2)t-\frac{(\lambda-a_1t)(\lambda-a_2t)}{\mu}. \nonumber
\end{align}
%fff:=q^2*k1*x^3-(q*(q+1)*k1*la-q^2*k1*k2*mu*(la-a3)/la+q*(t1+t2)*t/la-(la-a1*t)*(la-a2*t)/mu/la)*x^2+(q*k1*la^2-q*k1*k2*mu*(la-a3)+(q+1)*(t1+t2)*t-(la-a1*t)*(la-a2*t)/mu)*x-la*(t1+t2)*t;
%factor(eval(fff,la=a1*t));
%(-a1*t+x*q)*(k1*a1*t*q*x^2+x*q*k1*k2*mu*a1*t-x*q*k1*k2*mu*a3-x*q*k1*a1^2*t^2-x*t*t1-x*t*t2+a1*t^2*t1+a1*t^2*t2)/(a1*t)
%factor(eval(fff,la=a3));
%(-a3+x)*(q^2*k1*mu*a3*x^2+x*a3^2-x*a3*a2*t-x*q*t*mu*t2-x*q*k1*a3^2*mu-x*q*t*mu*t1-x*a1*t*a3+x*a1*t^2*a2+a3*t*mu*t2+a3*t*mu*t1)/(mu*a3)
We restrict the parameters to obtain the degenerated $q$-Heun equation.
We set $\lambda = a_3$.
Then, we have
\begin{align}
&  y_1(q^2x)-\{q^2\kappa_1x^2-dx+(\theta_1+\theta_2)t\}y_1(qx)
     +\kappa_1\kappa_2(qx-a_3)(x-a_1t)(x-a_2t)y_1(x)=0, \\
& d=q\kappa_1a_3+\frac{q\left(\theta_1+\theta_2\right)t}{a_3}-\frac{(a_3-a_1t)(a_3-a_2t)}{\mu a_3}. \nonumber
\end{align}
Note that the parameter $\mu $ appears only in $d$.
Let $u(x)$ be a function which satisfies $u(qx)=(x-a_1t)u(x)$ and set $f(x)=y_1(qx)/u(qx)$.
Then, the function $f(x)$ satisfies
\begin{align}
&     (qx-a_1t)f(qx)-\{q^2\kappa_1x^2+dx+(\theta_1+\theta_2)t\}f(x)
    +\kappa_1\kappa_2(qx-a_3)(x-a_2t)f(x/q)=0.
\end{align}

\subsection{$q$-P($A_5$)}
On the equation $Y(qx,t)=A(x,t)Y(x,t)$ for $q$-P($A_5$) by Murata, the matrix $A(x,t) $ is defined as 
\begin{align}
& A(x,t)=\begin{pmatrix}\kappa_1\{(x-\lambda)(x-\alpha)+\mu_1\}&\omega(x-\lambda)\\
          \kappa_1\omega^{-1}(\gamma x+\delta)&\kappa_2(x-\lambda)+\mu_2\end{pmatrix},\\
&  \mu_1=\frac{(\lambda-a_1t)(\lambda-a_2t)}{q\kappa_1\mu}, \; \mu_2=q\kappa_1\kappa_2\mu\lambda, \; \alpha=\frac{\lambda^{-1}\left\{\theta_1t-\kappa_1\mu_1-\mu_2\right\}+\kappa_2}{\kappa_1}, \nonumber \\
 & \gamma=\mu_2-\kappa_2\left\{(2\lambda+\alpha)-(a_1+a_2)t\right\},\; \delta=\lambda^{-1}(\alpha\lambda+\mu_1)(\kappa_2\lambda-\mu_2). \nonumber
\end{align}
Note that it follows that $\det A(x,t)=\kappa_1\kappa_2 x (x-a_1t)(x-a_2t)$ and the eigenvalues of the matrix $A(0,t)$ are $\theta_1t, 0 $.

Set $Y(x,t)= \! ^t (y_1(x), y_2(x))$.
It follows from the equation $Y(qx,t)=A(x,t)Y(x,t)$ that the function $y_1(x) $ satisfies the single $q$-difference equation
\begin{align}
& y_1(q^2x)
  -\frac{q^2\kappa_1x^3-b_2x^2+b_1x-\lambda\theta_1t}{x-\lambda}y_1(qx) \label{eq:y1q2A5} \\
&  +\frac{\kappa_1\kappa_2x(qx-\lambda)(x-a_1t)(x-a_2t)}{x-\lambda}y_1(x)=0, \nonumber
\end{align}
where
\begin{align}
&  b_2 =q(q+1)\kappa_1\lambda -q^2\kappa_1\kappa_2\mu +\frac{q\theta_1t}{\lambda} -\frac{(\lambda-a_1t)(\lambda-a_2t)}{\mu},\\
&  b_1 =q\kappa_1\lambda^2-q\kappa_1\kappa_2\mu\lambda +(q+1)\theta_1t-\frac{(\lambda-a_1t)(\lambda-a_2t)}{\mu}. \nonumber
\end{align}
We restrict the parameters to obtain the degenerated $q$-Heun equation.
We set $\lambda = a_1t $.
Then, we have
\begin{align}
&  y_1(q^2x)-\frac{(qx-a_1t)\{q\kappa_1x^2-dx+\theta_1t\}}{x-a_1t}y_1(qx) +\kappa_1\kappa_2x(qx-a_1t)(x-a_2t)y_1(x)=0,\\
&  d=q\kappa_1a_1t+\frac{\theta_1}{a_1}-{q\kappa_1\kappa_2\mu}. \nonumber
\end{align}
Note that the parameter $\mu $ appears only in $d$.
Set $y_1(x)=(x/q-a_1t)y(x)$.
We have
\begin{align}
& y(q^2x)-(q\kappa_1x^2-dx+\theta_1t)y(qx) +\kappa_1\kappa_2x\left(\frac{x}{q}-a_1t\right)(x-a_2t)y(x)=0.
\end{align}
Let $u(x)$ be a function which satisfies $u(qx)=(x/q-a_1t)u(x)$ and set $f(x)=y_1(qx)/u(qx)$.
Then, the function $f(x)$ satisfies
\begin{align}
&  (x-a_1t)f(qx)-(q\kappa_1x^2-dx+\theta_1t)f(x) +\kappa_1\kappa_2x(x-a_2t)f(x/q)=0.
\end{align}

We comment on another spacialization.
We set 
\begin{equation}
\mu = \frac{a_1 a_2 t }{q \theta _1 } +d \lambda 
\end{equation}
in equation (\ref{eq:y1q2A5}).
%mu=a1*a2*t/q/th1+zz*la, la=0
As $\lambda \to 0$, we have
\begin{align}
&  y_1(q^2x)- \Big( q^2\kappa_1x^2- q \Big( \theta _1 t d + \theta _1 \Big(\frac{1}{a_1}+\frac{1}{a_2} \Big)-\frac{ a_1 a_2 \kappa _1 \kappa _2 t }{\theta _1} \Big) x +\theta_1t \Big) y_1(qx) \\
& +q \kappa_1\kappa_2 x (x-a_1 t)( x-a_2 t) y_1(x)=0. \nonumber
\end{align}
%mu=a1*a2*t/q/th1+zz*la
%factor(coeff(pppp,x,2)); expand(coeff(pppp,x,1)); factor(coeff(pppp,x,0));
%q^2*k1
%a1*a2*q*k1*k2*t/th1-th1*q/a2-th1*t*q*zz-th1*q/a1
%th1*t

\subsection{$q$-P$(A_5)^\#$}
On the equation $Y(qx,t)=A(x,t)Y(x,t)$ for $q$-P$(A_5)^\#$ by Murata, the matrix $A(x,t) $ is defined as 
\begin{align}
& A(x,t)=\begin{pmatrix}\kappa_1\{(x-\lambda)(x-\alpha)+\mu_1\}&\omega(x-\lambda)\\
          \kappa_1\omega^{-1}(\gamma x+\delta)&\kappa_2(x-\lambda)+\mu_2\end{pmatrix},\\
 &  \mu_1 = \frac{\lambda(\lambda-a_1t)}{q\kappa_1\mu}, \; \mu_2 = q\kappa_1\kappa_2\mu(\lambda-a_3),\; \alpha=\frac{\lambda^{-1}\left(\theta_1t-\kappa_1\mu_1-\mu_2\right)+\kappa_2}{\kappa_1}, \nonumber \\
&  \gamma=\mu_2-\kappa_2\left(2\lambda+\alpha-a_1t-a_3\right), \; \delta=\lambda^{-1}(\alpha\lambda+\mu_1)(\kappa_2\lambda-\mu_2). \nonumber
\end{align}
Note that it follows that $\det A(x,t)=\kappa_1\kappa_2 x(x-a_1t)(x-a_3)$ and the eigenvalues of the matrix $A(0,t)$ are $\theta_1t, 0$.

Set $Y(x,t)= \! ^t (y_1(x), y_2(x))$.
It follows from the equation $Y(qx,t)=A(x,t)Y(x,t)$ that the function $y_1(x) $ satisfies the single $q$-difference equation
\begin{align}
& y_1(q^2x)
  -\frac{q^2\kappa_1x^3-b_2x^2+b_1x-\lambda\theta_1t}{x-\lambda}y_1(qx) \\
&  +\frac{\kappa_1\kappa_2x(qx-\lambda)(x-a_1t)(x-a_3)}{x-\lambda}y_1(x)=0,\nonumber 
\end{align}
where
\begin{align}
& b_2= q(q+1)\kappa_1\lambda-\frac{q^2\kappa_1\kappa_2\mu(\lambda-a_3)}{\lambda}  +\frac{q\theta_1t}{\lambda}-\frac{\lambda-a_1t}{\mu},\\
&  b_1= q\kappa_1\lambda^2-q\kappa_1\kappa_2\mu(\lambda-a_3)+(q+1)\theta_1t -\frac{\lambda(\lambda-a_1t)}{\mu}. \nonumber
\end{align}
We restrict the parameters to obtain the degenerated $q$-Heun equation.
We set $\lambda = a_3$.
Then, we have
\begin{align}
&  y_1(q^2x)-(q^2\kappa_1x^2-dx+\theta_1t)y_1(qx) +\kappa_1\kappa_2x(qx-a_3)(x-a_1t)y_1(x)=0,\\
&  d=q\kappa_1a_3+\frac{q\theta_1t}{a_3}-\frac{a_3-a_1t}{\mu}. \nonumber
\end{align}
Note that the parameter $\mu $ appears only in $d$.
Let $u(x)$ be a function which satisfies $u(qx)=(x-a_1t)u(x)$ and set $f(x)=y_1(qx)/u(qx)$.
Then, the function $f(x)$ satisfies
\begin{align}
&  (qx-a_1t)f(qx)-(q^2\kappa_1x^2+dx+\theta_1t)f(x) +\kappa_1\kappa_2x(qx-a_3)f(x/q)=0.
\end{align}

\subsection{$q$-P($A_6$)}
On the equation $Y(qx,t)=A(x,t)Y(x,t)$ for $q$-P($A_6$) by Murata, the matrix $A(x,t) $ is defined as 
\begin{align}
& A(x,t)=\begin{pmatrix}\kappa_1\{(x-\lambda)(x-\alpha)+\mu_1\}&\omega(x-\lambda)\\
          \kappa_1\omega^{-1}(\gamma x+\delta)&\kappa_2(x-\lambda)+\mu_2\end{pmatrix},\\
&  \mu_1=\frac{\lambda(\lambda-a_1t)}{q\kappa_1\mu}, \;  \mu_2=q\kappa_1\kappa_2\mu\lambda,\; \alpha=\frac{\lambda^{-1}\left\{\theta_1t-\kappa_1\mu_1-\mu_2\right\}+\kappa_2}{\kappa_1}, \nonumber \\
&  \gamma=\mu_2-\kappa_2(2\lambda+\alpha-a_1t), \; \delta=\lambda^{-1}(\alpha\lambda+\mu_1)(\kappa_2\lambda-\mu_2). \nonumber
\end{align}
Note that it follows that $\det A(x,t)=\kappa_1\kappa_2 x^2 (x-a_1t)$ and the eigenvalues of the matrix $A(0,t)$ are $\theta_1t, 0$.

Set $Y(x,t)= \! ^t (y_1(x), y_2(x))$.
It follows from the equation $Y(qx,t)=A(x,t)Y(x,t)$ that the function $y_1(x) $ satisfies the single $q$-difference equation
\begin{align}
& y_1(q^2x) -\frac{q^2\kappa_1x^3-b_2x^2+b_1x-\lambda\theta_1t}{x-\lambda}y_1(qx) +\frac{\kappa_1\kappa_2x^2(qx-\lambda)(x-a_1t)}{x-\lambda}y_1(x)=0, \label{eq:y1q2A6}
\end{align}
where
\begin{align}
& b_2=q(q+1)\kappa_1\lambda-q^2\kappa_1\kappa_2\mu
      +\frac{q\theta_1t}{\lambda}-\frac{\lambda-a_1t}{\mu},\\
&  b_1=q\kappa_1\lambda^2-q\kappa_1\kappa_2\mu\lambda+(q+1)\theta_1t
      -\frac{\lambda(\lambda-a_1t)}{\mu}. \nonumber
\end{align}
We restrict the parameters to obtain the degenerated $q$-Heun equation.
We set $\lambda = a_1t $.
Then, we have
  \begin{align}
&  y_1(q^2x)-\frac{(qx-a_1t)(q\kappa_1x^2-dx+\theta_1t)}{x-a_1t}y_1(qx) +\kappa_1\kappa_2x^2(qx-a_1t)y_1(x)=0,\\
&  d=q\kappa_1a_1t+\frac{\theta_1}{a_1}-q\kappa_1\kappa_2\mu. \nonumber
  \end{align}
Note that the parameter $\mu $ appears only in $d$.
Set $y_1(x)=(x/q-a_1t)y(x)$.
We have
\begin{align}
  y(q^2x)-(q\kappa_1x^2-dx+\theta_1t)y(qx)
  +\kappa_1\kappa_2x^2\left(\frac{x}{q}-a_1t\right)y(x)=0.
\end{align}
Let $u(x)$ be a function which satisfies $u(qx)=(x/q-a_1t)u(x)$ and set $f(x)=y_1(qx)/u(qx)$.
Then, the function $f(x)$ satisfies
  \begin{align}
    (x-a_1t)f(qx)-(q\kappa_1x^2-dx+\theta_1t)f(x)
    +\kappa_1\kappa_2x^2f(x/q)=0.
  \end{align}

We comment on another spacialization.
We set 
\begin{equation}
\mu = \frac{\lambda (\lambda -a_1 t) }{q \theta _1 t}(1 +d \lambda ) 
\end{equation}
in equation (\ref{eq:y1q2A6}).
%mu=la*(la-a1*t)/q/th1/t*(1+zz*la)), la=0
As $\lambda \to 0$, we have
 \begin{align}
    y_1(q^2x)- (q^2\kappa_1x^2- q \theta _1 t d  x +\theta_1t)y_1(qx) +q \kappa_1\kappa_2 x^2 (x-a_1 t) y_1(x)=0 .
  \end{align}

\subsection{$q$-P$(A_6)^\#$}
On the equation $Y(qx,t)=A(x,t)Y(x,t)$ for $q$-P$(A_6)^\#$ by Murata, the matrix $A(x,t) $ is defined as 
\begin{align}
&    A(x,t)=\begin{pmatrix}\kappa_1\{(x-\lambda)(x-\alpha)+\mu_1\}&\omega(x-\lambda)\\
          \kappa_1\omega^{-1}(\gamma x+\delta)&\kappa_2(x-\lambda)+\mu_2\end{pmatrix},\\
&    \mu_1=\frac{\lambda ^2}{q\kappa_1\mu}, \;
    \mu_2=q\kappa_1\kappa_2\mu (\lambda -a_3) ,\;
    \alpha=\frac{\lambda^{-1}\left\{\theta_1t-\kappa_1\mu_1-\mu_2\right\}+\kappa_2}{\kappa_1}, \nonumber \\
 &   \gamma=\mu_2-\kappa_2(2\lambda+\alpha-a_3),
    \delta=\lambda^{-1}(\alpha\lambda+\mu_1)(\kappa_2\lambda-\mu_2). \nonumber
\end{align}
Note that it follows that $\det A(x,t)=\kappa_1\kappa_2x^2(x-a_3)$ and the eigenvalues of the matrix $A(0,t)$ are $\theta_1t,0$.

Set $Y(x,t)= \! ^t (y_1(x), y_2(x))$.
It follows from the equation $Y(qx,t)=A(x,t)Y(x,t)$ that the function $y_1(x) $ satisfies the single $q$-difference equation
\begin{align}
  y_1(q^2x)
  -\frac{q^2\kappa_1x^3-b_2x^2+b_1x-\lambda\theta_1t}{x-\lambda}y_1(qx)
  +\frac{\kappa_1\kappa_2x^2(qx-\lambda)(x-a_3)}{x-\lambda}y_1(x)=0,
 \end{align}
where
\begin{align}
&  b_2 =q(q+1)\kappa_1\lambda-\frac{q^2\kappa_1\kappa_2\mu(\lambda-a_3)}{\lambda}
        +\frac{q\theta_1t}{\lambda}-\frac{\lambda}{\mu},\\
&  b_1 =q\kappa_1\lambda^2-q\kappa_1\kappa_2\mu(\lambda-a_3)+(q+1)\theta_1t
        -\frac{\lambda^2}{\mu}. \nonumber 
\end{align}
We restrict the parameters to obtain the degenerated $q$-Heun equation.
We set $\lambda = a_3$.
Then, we have
 \begin{align}
 & y_1(q^2x)-(q^2\kappa_1x^2-dx+\theta_1t)y_1(qx)+\kappa_1\kappa_2x^2(qx-a_3)y_1(x)=0, \\
 &  d=q\kappa_1a_3+\frac{q\theta_1t}{a_3}-\frac{a_3}{\mu}. \nonumber
  \end{align}
Note that the parameter $\mu $ appears only in $d$.
Let $u(x)$ be a function which satisfies $u(qx)=(qx-a_3)u(x)$ and set $f(x)=y_1(qx)/u(qx)$.
Then, the function $f(x)$ satisfies
  \begin{align}
    (q^2x-a_3)f(qx)-(q^2\kappa_1x^2-dx+\theta_1t)f(x)+\kappa_1\kappa_2x^2f(x/q)=0.
  \end{align}

\subsection{$q$-P($A_7$)}
On the equation $Y(qx,t)=A(x,t)Y(x,t)$ for $q$-P($A_7$) by Murata, the matrix $A(x,t) $ is defined as 
\begin{align}
%}\renewcommand{\arraystretch}{1.2}
 &   A(x,t)=\begin{pmatrix}\kappa_1\{(x-\lambda)(x-\alpha)+\mu_1\}&\omega(x-\lambda)\\
          \kappa_1\omega^{-1}(\gamma x+\delta)&\kappa_2(x-\lambda)+\mu_2\end{pmatrix},\\
&   \mu_1= \frac{\lambda^2}{q\kappa_1\mu}, \; \mu_2 = q\kappa_1\kappa_2\mu\lambda, \; \alpha=\frac{\lambda^{-1}(\theta_1t-\kappa_1\mu_1-\mu_2)+\kappa_2}{\kappa_1}, \nonumber \\
&    \gamma=\mu_2-\kappa_2(2\lambda+\alpha), \;  \delta=\lambda^{-1}(\alpha\lambda+\mu_1)(\kappa_2\lambda-\mu_2). \nonumber
\end{align}
Note that it follows that $\det A(x,t)=\kappa_1\kappa_2x^3$ and the eigenvalues of the matrix $A(0,t)$ are $\theta_1t,0$.

Set $Y(x,t)= \! ^t (y_1(x), y_2(x))$.
It follows from the equation $Y(qx,t)=A(x,t)Y(x,t)$ that the function $y_1(x) $ satisfies the single $q$-difference equation
\begin{align}
  y_1(q^2x)
  -\frac{q^2\kappa_1x^3-b_2x^2+b_1x-\lambda\theta_1t}{x-\lambda}y_1(qx)
  +\frac{\kappa_1\kappa_2x^3(qx-\lambda)}{x-\lambda}y_1(x)=0,
\end{align}
where
\begin{align}
&   b_2 =q(q+1)\kappa_1\lambda-q^2\kappa_1\kappa_2\mu +\frac{q\theta_1t}{\lambda}-\frac{\lambda}{\mu},\\
&  b_1 =q\kappa_1\lambda^2-q\kappa_1\kappa_2\mu \lambda +(q+1)\theta_1t
        -\frac{\lambda^2}{\mu}. \nonumber 
\end{align}
We restrict the parameters to obtain the degenerated $q$-Heun equation.
We set 
\begin{equation}
\mu = \frac{\lambda^2 }{q \theta _1 t} (1+d \lambda ) .
\end{equation}
%mu=la^2/q/th1/t*(1+zz*la), la=0
As $\lambda \to 0$, we have
 \begin{align}
    y_1(q^2x)-(q^2\kappa_1x^2- q \theta _1 t d x+\theta_1t)y_1(qx)+q \kappa_1\kappa_2x^3 y_1(x)=0 .
  \end{align}
Let $u(x)$ be a function which satisfies $u(qx)=x u(x)$ and set $f(x)=y_1(qx)/u(qx)$.
Then, the function $f(x)$ satisfies
\begin{align}
  qxf(qx)-(q^2\kappa_1x^2- q \theta _1 t d x+\theta_1t)f(x)+q\kappa_1\kappa_2x^2f(x/q)=0.
\end{align}

\subsection{$q$-P($A_7'$)}
On the equation $Y(qx,t)=A(x,t)Y(x,t)$ for $q$-P($A_7'$) by Murata, the matrix $A(x,t) $ is defined as 
\begin{align}
%}\renewcommand{\arraystretch}{1.2}
&    A(x,t)=\begin{pmatrix}\kappa_1\{(x-\lambda)(x-\alpha)+\mu_1\}&\omega(x-\lambda)\\
          \kappa_1\omega^{-1}(\gamma x+\delta)&\mu_2\end{pmatrix},\\
 &   \mu_1 = \frac{\lambda^2}{q\kappa_1\mu}, \;  \mu_2 = q\kappa_1\kappa_2\mu, \;  \alpha=\frac{\lambda^{-1}(\theta_1t-\kappa_1\mu_1-\mu_2)}{\kappa_1}, \nonumber \\
&   \gamma=\mu_2 - \kappa_2, \;  \delta=-\lambda^{-1}\mu_2(\alpha\lambda+\mu_1). \nonumber
\end{align}
Note that it follows that $\det A(x,t)=\kappa_1\kappa_2x^2$ and the eigenvalues of the matrix $A(0,t)$ are $\theta_1t,0$.

Set $Y(x,t)= \! ^t (y_1(x), y_2(x))$.
It follows from the equation $Y(qx,t)=A(x,t)Y(x,t)$ that the function $y_1(x) $ satisfies the single $q$-difference equation
\begin{align}
  y_1(q^2x)
  -\frac{q^2\kappa_1x^3-b_2x^2+b_1x-\lambda\theta_1t}{x-\lambda}y_1(qx)
  +\frac{\kappa_1\kappa_2x^2(qx-\lambda)}{x-\lambda}y_1(x)=0,
\end{align}
where
\begin{align}
&  b_2 =q(q+1)\kappa_1\lambda-\frac{q^2\kappa_1\kappa_2\mu }{\lambda } +\frac{q\theta_1t}{\lambda}-\frac{\lambda}{\mu},\\
&  b_1 =q\kappa_1\lambda^2-q\kappa_1\kappa_2\mu +(q+1)\theta_1t
        -\frac{\lambda^2}{\mu}. \nonumber 
\end{align}
We restrict the parameters to obtain the degenerated $q$-Heun equation.
We set 
\begin{equation}
\mu = \frac{\theta _1 t }{q \kappa_1\kappa_2 } +d \lambda  .
\end{equation}
%mu=th1*t/q/k1/k2+la*zz, la=0
As $\lambda \to 0$, we have
 \begin{align}
    y_1(q^2x)-q (q \kappa_1x^2 + q \kappa_1\kappa_2  d x+\theta_1t)y_1(qx)+q \kappa_1\kappa_2x^2 y_1(x)=0 .
  \end{align}
Set $f(x)=y_1(qx)$.
Then, the function $f (x)$ satisfies
\begin{align}
    f(qx)-q (q \kappa_1x^2 + q \kappa_1\kappa_2  d x+\theta_1t)f(x)+q \kappa_1\kappa_2x^2 f(x/q)=0 .
\end{align}

\subsection{} \label{sec:Mlist}

We summarize the degenerated $q$-Heun equation obtained from the list of Murata.
%\vspace{0.2cm}
\begin{center}  
    \begingroup\renewcommand{\arraystretch}{1.2}
      \begin{tabular}{lll}
        $q$-P$(A_4)$
        & $(qx-a_1t)f(qx)-\{q^2\kappa_1x^2+dx+(\theta_1+\theta_2)t\}f(x) $ \\
& $ \qquad \qquad \qquad +\kappa_1\kappa_2(qx-a_3)(x-a_2t)f(x/q)=0 , $\\
%        $(q^2x-a_1t)f(qx)-\{q^2\kappa_1x^2-dx+(\theta_1+\theta_2)t\}f(x)   +\kappa_1\kappa_2(x-a_2t)(x-a_3)f(x/q)=0$\\
        $q$-P$(A_5)$
        & $(x-a_1t)f(qx)-(q\kappa_1x^2-dx+\theta_1t)f(x) $ \\
& $ \qquad \qquad \qquad \qquad \qquad +\kappa_1\kappa_2x(x-a_2t)f(x/q)=0 , $\\
%        &$(x-a_1t)f(qx)-(q\kappa_1x^2-dx+\theta_1t)f(x) +\kappa_1\kappa_2x(x-a_2t)f(x/q)=0$\\
        $q$-P$(A_5)^\#$
        &$(qx-a_1t)f(qx)-(q^2\kappa_1x^2+dx+\theta_1t)f(x)  $ \\
& $ \qquad \qquad \qquad \qquad \qquad +\kappa_1\kappa_2x(qx-a_3)f(x/q)=0 , $\\
%        &$(x-a_1t)f(qx)-(q\kappa_1x^2-dx+\theta_1t)f(x) +\kappa_1\kappa_2x(x-a_3)f(x/q)=0$\\
        $q$-P$(A_6)$
        &$(x-a_1t)f(qx)-(q\kappa_1x^2-dx+\theta_1t)f(x)
        +\kappa_1\kappa_2x^2f(x/q)=0 , $\\
        $q$-P$(A_6)^\#$
        &$(q^2x-a_3)f(qx)-(q^2\kappa_1x^2-dx+\theta_1t)f(x)+\kappa_1\kappa_2x^2f(x/q)=0 , $\\
        $q$-P$(A_7)$
        &$qxf(qx)-(q^2\kappa_1x^2- q \theta _1 t d x+\theta_1t)f(x)+q\kappa_1\kappa_2x^2f(x/q)=0 , $\\
%$qxf(qx)-(q^2\kappa_1x^2+qdx+\theta_1t)f(x)+q\kappa_1\kappa_2x^2f(x/q)=0$\\
        $q$-P$(A_7')$
        &$f(qx)-q (q \kappa_1x^2 + q \kappa_1\kappa_2  d x+\theta_1t)f(x)+q \kappa_1\kappa_2x^2 f(x/q)=0 . $
%qf(qx)-(q\kappa_1x^2-qdx+q^2\theta_1t)f(x)+\kappa_1\kappa_2x^2f(x/q)=0$
      \end{tabular}
    \endgroup
%\vspace{0.2cm}
\end{center}
The Newton diagrams of these equations are described as follows.
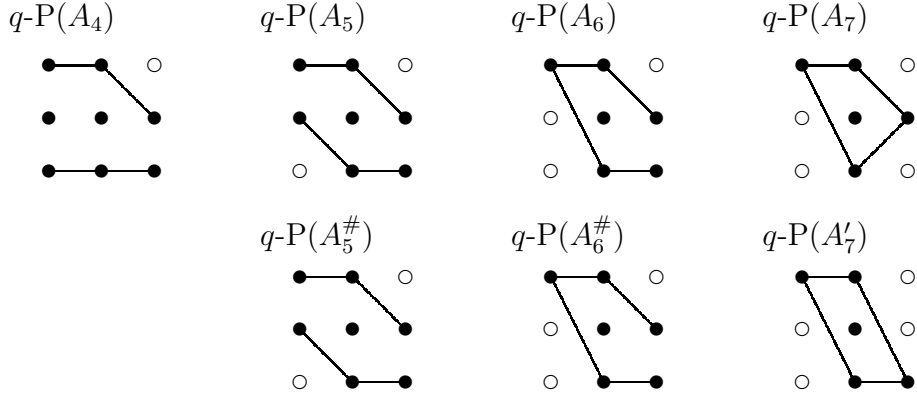
\begin{center}  
\begin{picture}(370,175)(0,0)
  \put(15,155){$q$-P($A_4$)}
    \put(30,140){\circle*{5}}
    \put(50,140){\circle*{5}}
    \put(70,140){\circle{5}}
    \put(30,120){\circle*{5}}
    \put(50,120){\circle*{5}}
    \put(70,120){\circle*{5}}
    \put(30,100){\circle*{5}}
    \put(50,100){\circle*{5}}
    \put(70,100){\circle*{5}}
%    \qbezier(30,140)(30,120)(30,100)
    \qbezier(30,100)(50,100)(70,100)
%    \qbezier(70,100)(70,110)(70,120)
    \qbezier(70,120)(60,130)(50,140)
    \qbezier(50,140)(40,140)(30,140)
  \put(110,155){$q$-P($A_5$)}
    %\put(90,117.5){\vector(1,0){15}}
    \put(125,140){\circle*{5}}
    \put(145,140){\circle*{5}}
    \put(165,140){\circle{5}}
    \put(125,120){\circle*{5}}
    \put(145,120){\circle*{5}}
    \put(165,120){\circle*{5}}
    \put(125,100){\circle{5}}
    \put(145,100){\circle*{5}}
    \put(165,100){\circle*{5}}
    \qbezier(125,140)(135,140)(145,140)
    \qbezier(145,140)(155,130)(165,120)
%    \qbezier(165,120)(165,110)(165,100)
    \qbezier(165,100)(155,100)(145,100)
    \qbezier(145,100)(135,110)(125,120)
%    \qbezier(125,120)(125,130)(125,140)
  \put(110,72){$q$-P($A_5^\#$)}
    %\put(90,110){\vector(1,-4){15}}
    \put(125,60){\circle*{5}}
    \put(145,60){\circle*{5}}
    \put(165,60){\circle{5}}
    \put(125,40){\circle*{5}}
    \put(145,40){\circle*{5}}
    \put(165,40){\circle*{5}}
    \put(125,20){\circle{5}}
    \put(145,20){\circle*{5}}
    \put(165,20){\circle*{5}}
    \qbezier(125,60)(135,60)(145,60)
    \qbezier(145,60)(155,50)(165,40)
%    \qbezier(165,40)(165,30)(165,20)
    \qbezier(165,20)(155,20)(145,20)
    \qbezier(145,20)(135,30)(125,40)
%    \qbezier(125,40)(125,50)(125,60)
  \put(205,155){$q$-P($A_6$)}
    %\put(185,117.5){\vector(1,0){15}}
    \put(220,140){\circle*{5}}
    \put(240,140){\circle*{5}}
    \put(260,140){\circle{5}}
    \put(220,120){\circle{5}}
    \put(240,120){\circle*{5}}
    \put(260,120){\circle*{5}}
    \put(220,100){\circle{5}}
    \put(240,100){\circle*{5}}
    \put(260,100){\circle*{5}}
    \qbezier(220,140)(230,140)(240,140)
    \qbezier(240,140)(250,130)(260,120)
%    \qbezier(260,120)(260,110)(260,100)
    \qbezier(260,100)(250,100)(240,100)
    \qbezier(220,140)(230,120)(240,100)
  \put(205,72){$q$-P($A_6^\#$)}
    %\put(185,45){\vector(1,4){15}}
    %\put(185,37.5){\vector(1,0){15}}
    \put(220,60){\circle*{5}}
    \put(240,60){\circle*{5}}
    \put(260,60){\circle{5}}
    \put(220,40){\circle{5}}
    \put(240,40){\circle*{5}}
    \put(260,40){\circle*{5}}
    \put(220,20){\circle{5}}
    \put(240,20){\circle*{5}}
    \put(260,20){\circle*{5}}
    \qbezier(220,60)(230,60)(240,60)
    \qbezier(240,60)(250,50)(260,40)
%    \qbezier(260,40)(260,30)(260,20)
    \qbezier(260,20)(250,20)(240,20)
    \qbezier(220,60)(230,40)(240,20)
  \put(300,155){$q$-P($A_7$)}
    %\put(280,117.5){\vector(1,0){15}}
    \put(315,140){\circle*{5}}
    \put(335,140){\circle*{5}}
    \put(355,140){\circle{5}}
    \put(315,120){\circle{5}}
    \put(335,120){\circle*{5}}
    \put(355,120){\circle*{5}}
    \put(315,100){\circle{5}}
    \put(335,100){\circle*{5}}
    \put(355,100){\circle{5}}
    \qbezier(335,140)(345,130)(355,120)
    \qbezier(355,120)(345,110)(335,100)
    \qbezier(315,140)(325,120)(335,100)
    \qbezier(315,140)(325,140)(335,140)
  \put(300,72){$q$-P($A_7'$)}
    %\put(280,45){\vector(1,4){15}}
    %\put(280,37.5){\vector(1,0){15}}
    \put(315,60){\circle*{5}}
    \put(335,60){\circle*{5}}
    \put(355,60){\circle{5}}
    \put(315,40){\circle{5}}
    \put(335,40){\circle*{5}}
    \put(355,40){\circle{5}}
    \put(315,20){\circle{5}}
    \put(335,20){\circle*{5}}
    \put(355,20){\circle*{5}}
    \qbezier(315,60)(325,60)(335,60)
    \qbezier(355,20)(345,20)(335,20)
    \qbezier(355,20)(345,40)(335,60)
    \qbezier(335,20)(325,40)(315,60)
  \end{picture}
\end{center}  
In section \ref{sec:cqHs}, we introduce the confluent $q$-Heun equation, the biconfluent $q$-Heun equation and the doubly confluent $q$-Heun equation.
The equation corresponding to $q$-P$(A_4)$ is the confluent $q$-Heun equation, and the equations corresponding to $q$-P$(A_5)$, $q$-P$(A_5)^\#$ are the doubly confluent $q$-Heun equation.

\section{Degenerations of $q$-Heun equations arising from the list of Kajiwara-Noumi-Yamada} \label{sec:KNY}

In \cite{KNY}, Kajiwara, Noumi and Yamada presented a list of Lax pairs $L_1$ and $L_2$ for discrete Painlev\'e equations from which the discrete Painlev\'e equations are derived by the compatibility condition.

In this section, we obtain degenerate $q$-Heun equations from the linear equation $L_1y(z,t)=0$.
The non-degenerate $q$-Heun equation was obtained from the equation $L_1y(z,t)=0$ for $q$-P($D_5^{(1)}$) in \cite{TakR}, and we discuss the case $q$-P($A_4^{(1)}$) and its descendants.
In the following, we assume that the parameters $\kappa_1, \kappa_2, \nu_1, \nu _2, \dots , \nu _8 $ satisfy $ \kappa_1^2\kappa_2^2 = q  \nu_1 \nu _2 \cdots \nu _8 $.

\subsection{$q$-P($A_4^{(1)}$)}

Let $T_z$ be the $q$-shift operator on the variable $z$ such that $T_z y(z)=y(qz)$.
On the equation $L_1y(z,t)=0$ for $q$-P($A_4^{(1)}$) by Kajiwara, Noumi, Yamada, the operator $L_1 $ is defined as 
\begin{align}
  L_1=& \frac{\nu_1 \nu_2 \nu_3 \nu_4 (g-\frac{\nu_5}{\kappa_2}) (g-\frac{\nu_6}{\kappa_2})}{fg}
        +\frac{(g\nu_1-1)z}{qg} \\
&  +\frac{ \nu_1 \nu_2 \nu_3 (\frac zq-\nu_4 )}{f-\frac zq}(g-T_z^{-1})
        +\frac{(z-\frac{\kappa_1}{\nu_7}) (z-\frac{\kappa_1}{\nu_8}) }{q(f-z)}\left(T_z-\frac1g\right). \nonumber
%    q\prod_{i=1}^8\nu_i=\kappa_1^2\kappa_2^2.
\end{align}
We set $f=\nu_4$.
Then, the equation $L_1 y(z)= 0$ is written as
\begin{align}
&  {\left(z-\frac{\kappa_1}{\nu_7}\right) \left(z-\frac{\kappa_1}{\nu_8}\right) }y(qz)
   +\left(-{\nu_1}z^2+dz-\frac{\kappa_1^2\kappa_2(\nu_5+\nu_6)}{\nu_5 \nu_6 \nu_7 \nu_8 }\right)y(z) \\
&    +q\nu_1 \nu_2 \nu_3 (z-\nu_4)y(z/q)=0, \nonumber \\
&   d=\frac{\nu_1\left\{q\nu_2\nu_3(\nu_5+\nu_6)+\kappa_2\nu_4\right\}}{q\kappa_2}
      +\frac{\left(\kappa_1-\nu_4\nu_7\right)\left(\kappa_1-\nu_4\nu_8\right)}{q\nu_4\nu_7\nu_8}\frac1g. \nonumber 
\end{align}
%where the term $d$ is expressed by using the parameter $g$.

\subsection{$q$-P($E_3^{(1)};a$)}
On the equation $L_1y(z,t)=0$ for $q$-P($E_3^{(1)};a$) by Kajiwara, Noumi, Yamada, the operator $L_1 $ is defined as 
\begin{align}
  L_1= & \frac{\nu_1 \nu_2 \nu_3 \nu_4 (g-\frac{\nu_5}{\kappa_2}) (g-\frac{\nu_6}{\kappa_2}) }{fg}
        +\frac{\nu_1z}{q}
        +\frac{\nu_1 \nu_2 \nu_3 (\frac zq-\nu_4 )}{f-\frac zq}(g-T_z^{-1}) \\
&        -\frac{\frac{\kappa_1}{\nu_8} (z-\frac{\kappa_1}{\nu_7})}{q(f-z)}
          \left(T_z-\frac1g\right). \nonumber
\end{align}
We set $f=\nu_4$.
Then, the equation $L_1 y(z)= 0$ is written as
\begin{align}
&  \frac{\kappa_1(z-\frac{\kappa_1}{\nu_7})}{\nu_8}y(qz)
    +\left(\nu_1z^2+dz+\frac{\kappa_1^2\kappa_2(\nu_5+\nu_6)}{\nu_5 \nu_6 \nu_7 \nu_8 }\right)y(z)
    +q\nu_1 \nu_2 \nu_3 (z-\nu_4)y(z/q)=0, \\
&    d=-\frac{\nu_1\left\{\nu_2\nu_3(\nu_5+\nu_6)+\kappa_2\nu_4\right\}}{\kappa_2}
      +\frac{\kappa_1(\kappa_1-\nu_4\nu_7)}{q\nu_4\nu_7\nu_8}\frac1g. \nonumber 
\end{align}
%where the term $d$ is expressed by using the parameter $g$.
By a gauge-transformation, we obtain
\begin{align}
&  \frac{\kappa_1}{\nu_8} (qz-\nu_4 ) \left(z-\frac{\kappa_1}{\nu_7}\right) y(qz)
    +\left(\nu_1z^2+dz+\frac{\kappa_1^2\kappa_2(\nu_5+\nu_6)}{\nu_5 \nu_6 \nu_7 \nu_8 }\right)y(z)
    +q\nu_1 \nu_2 \nu_3 y(z/q)=0. 
\end{align}
%\begin{gather*}
%  \frac{\kappa_1}{\nu_8}(q^2z-\nu_4)\left(qz-\frac{\kappa_1}{\nu_7}\right)y(qz)
%  +\left(q^2\nu_1z^2+qdz+\frac{\kappa_1^2\kappa_2(\nu_5+\nu_6)}{\prod_{i=5}^8\nu_i}\right)y(z)
%  +q\left(\prod_{i=1}^3\nu_i\right)y(z/q)=0.
%\end{gather*}

\subsection{$q$-P($E_3^{(1)};b$)}
On the equation $L_1y(z,t)=0$ for $q$-P($E_3^{(1)};b$) by Kajiwara, Noumi, Yamada, the operator $L_1 $ is defined as 
\begin{align}
    L_1= & \left(g-\frac{\nu_5}{\kappa_2}\right)\frac{\nu_1 \nu_2 \nu_3 \nu_4 }{f}
        +\frac{(g\nu_1-1)z}{qg}
        +\frac{\nu_1 \nu_2 \nu_3 (\frac zq-\nu_4 )}{f-\frac zq}(g-T_z^{-1}) \\
&        +\frac{z(z-\frac{\kappa_1}{\nu_8})}{q(f-z)}\left(T_z-\frac1g\right) . \nonumber 
\end{align}
We set $f=\nu_4$.
Then, the equation $L_1 y(z)= 0$ is written as
\begin{align}
&   {z\left(z-\frac{\kappa_1}{\nu_8}\right)}y(qz)
    +\left(-{\nu_1}z^2+qdz-\frac{q\prod_{i=1}^5\nu_i}{\kappa_2}\right)y(z)
    +q\nu_1 \nu_2 \nu_3 (z-\nu_4)y(z/q)=0, \\
& d=\frac{\nu_1(q\nu_2\nu_3\nu_5+\kappa_2\nu_4)}{q\kappa_2} +\frac{\kappa_1-\nu_4\nu_8}{q\nu_8} \frac1g. \nonumber
\end{align}
%where the term $d$ is expressed by using the parameter $g$.

\subsection{$q$-P($E_2^{(1)};a$)}
On the equation $L_1y(z,t)=0$ for $q$-P($E_2^{(1)};a$) by Kajiwara, Noumi, Yamada, the operator $L_1 $ is defined as 
\begin{equation}
    L_1=\left(g-\frac{\nu_5}{\kappa_2}\right)\frac{\nu_1 \nu_2 \nu_3 \nu_4 }{f}
        +\frac{\nu_1z}{q}
        +\frac{\nu_1 \nu_2 \nu_3 (\frac{z}{q}-\nu_4 )}{f-\frac{z}{q}}
          (g-T_z^{-1})
        -\frac{\frac{\kappa_1}{\nu_8}\cdot z}{q(f-z)}\left(T_z-\frac{1}{g}\right) .
\end{equation}
We set $f=\nu_4$.
Then, the equation $L_1 y(z)= 0$ is written as
\begin{align}
&   \frac{\kappa_1z}{\nu_8}y(qz)
    +\left(\nu_1z^2+dz+\frac{q\prod_{i=1}^5\nu_i}{\kappa_2}\right)y(z)
    +q\nu_1 \nu_2 \nu_3 (z-\nu_4)y(z/q)=0, \\
& d=\frac{\nu_1(q\nu_2\nu_3\nu_5+\kappa_2\nu_4)}{q\kappa_2}+\frac{\kappa_1}{q\nu_8}\frac1g . \nonumber
\end{align}
%where the term $d$ is expressed by using the parameter $g$.
By a gauge-transformation, we obtain
\begin{align}
&    \frac{\kappa_1}{\nu_8} z (qz-\nu_4) y(qz)
    +\left(\nu_1z^2+dz+\frac{q\prod_{i=1}^5\nu_i}{\kappa_2}\right)y(z)
    +q\nu_1 \nu_2 \nu_3 y(z/q)=0 .
\end{align}
%\begin{gather*}
%  \frac{\kappa_1}{\nu_8}z\left(q^2z-{\nu_4}\right)y(qz)
%  +\left(q\nu_1z^2+dz+\frac{\prod_{i=1}^5\nu_i}{\kappa_2}\right)y(z)
%  +\nu_1 \nu_2 \nu_3 y(z/q)=0.
%\end{gather*}

\subsection{$q$-P($E_2^{(1)};b$)}
On the equation $L_1y(z,t)=0$ for $q$-P($E_2^{(1)};b$) by Kajiwara, Noumi, Yamada, the operator $L_1 $ is defined as 
\begin{equation}
    L_1=\frac{g\nu_1 \nu_2 \nu_3 \nu_4 }{f}
        +\frac{(g\nu_1-1)z}{qg}
        +\frac{\nu_1 \nu_2 \nu_3 (\frac zq-\nu_4)}{f-\frac zq}(g-T_z^{-1})
        +\frac{z (z-\frac{\kappa_1}{\nu_8})}{q(f-z)}\left(T_z-\frac1g\right) .
\end{equation}
We set $f=\nu_4$.
Then, the equation $L_1 y(z)= 0$ is written as
\begin{align}
&   z\left(z-\frac{\kappa_1}{\nu_8}\right)y(qz)
    +\left(-{\nu_1}z^2+dz\right)y(z)
    -q\nu_1 \nu_2 \nu_3 (z-\nu_4)y(z/q)=0, \\
& d=\frac{\nu_1\nu_4}{q}+\frac{\kappa_1-\nu_4\nu_8}{q\nu_8} \frac1g . \nonumber 
\end{align}
%where the term $d$ is expressed by using the parameter $g$.

\subsection{$q$-P($A_1^{(1)}$)}
On the equation $L_1y(z,t)=0$ for $q$-P($A_1^{(1)}$) by Kajiwara, Noumi, Yamada, the operator $L_1 $ is defined as 
\begin{equation}
    L_1=\frac{g\nu_1 \nu_2 \nu_3 \nu_4 }{f}
        -\frac{z}{qg}
        +\frac{\nu_1 \nu_2 \nu_3 (\frac zq-\nu_4)}{f-\frac zq}(g-T_z^{-1})
        +\frac{z (z-\frac{\kappa_1}{\nu_8})}{q(f-z)}\left(T_z-\frac1g\right).
\end{equation}
We set $f=\nu_4$.
Then, the equation $L_1 y(z)= 0$ is written as
\begin{align}
&    z\left(z-\frac{\kappa_1}{\nu_8}\right)y(qz)
    +dzy(z)
    -q\nu_1 \nu_2 \nu_3 (z-\nu_4)y(z/q)=0, \\ 
& d=\frac{\kappa_1-\nu_4\nu_8}{q\nu_8} \frac1g. \nonumber
\end{align}
%where the term $d$ is expressed by using the parameter $g$.

\subsection{Another $q$-P($A_1^{(1)}/A_7^{(1)}$)}
On the equation $L_1y(z,t)=0$ for $q$-P$\underset{|\alpha^2|=8}{(A_1^{(1)})}$ by Kajiwara, Noumi, Yamada, the operator $L_1 $ is defined as 
\begin{equation}
    L_1=\frac{g\nu_1 \nu_2 \nu_3 \nu_4 }{f}
        +\frac{\nu_1z}{q}
        +\frac{\nu_1 \nu_2 \nu_3 (\frac zq-\nu_4)}{f-\frac zq}(g-T_z^{-1})
        -\frac{\frac{\kappa_1}{\nu_8}\cdot z}{q(f-z)}\left(T_z-\frac1g\right).
\end{equation}
We set $f=\nu_4$.
Then, the equation $L_1 y(z)= 0$ is written as
\begin{align}
&    {\frac{\kappa_1}{\nu_8}}zy(qz)
    + (\nu_1z^2-dz ) y(z)
    +q\nu_1 \nu_2 \nu_3 (z-\nu_4)y(z/q)=0, \\
& d=\frac{\nu_1\nu_4}{q}+\frac{\kappa_1}{q\nu_8} \frac1g. \nonumber  
\end{align}
%where the term $d$ is expressed by using the parameter $g$.
By a gauge-transformation, we obtain
\begin{align}
&   {\frac{\kappa_1}{\nu_8}}z(qz-\nu_4)y(qz)
    + (\nu_1z^2-dz ) y(z)
    +q\nu_1\nu_2\nu_3y(z/q)=0.
\end{align}

\subsection{}

We summarize the degenerated $q$-Heun equation obtained from the list of Kajiwara-Noumi-Yamada.
%\vspace{0.2cm}
\begin{center}
    {\begingroup\renewcommand{\arraystretch}{2.3}
      \begin{tabular}{lll}
        $q$-P$(A_4^{(1)})$
           &$\displaystyle{
           \left(z-\frac{\kappa_1}{\nu_7}\right)\left(z-\frac{\kappa_1}{\nu_8}\right)y(qz)
             +\left(-{\nu_1}z^2+dz-\frac{\kappa_1^2\kappa_2(\nu_5+\nu_6)}{\nu_5 \nu_6 \nu_7 \nu_8 }\right)y(z)}$\\
          &$\displaystyle{
           +q\nu_1 \nu_2 \nu_3 (z-\nu_4)y(z/q)=0,}$\\
        $q$-P$(E_3^{(1)})(a)$
           &$\displaystyle
           {\frac{\kappa_1}{\nu_8} (qz-\nu_4 ) \left(z-\frac{\kappa_1}{\nu_7}\right) y(qz)
    +\left(\nu_1z^2+dz+\frac{\kappa_1^2\kappa_2(\nu_5+\nu_6)}{\nu_5 \nu_6 \nu_7 \nu_8 }\right)y(z)}$\\
           &$\displaystyle
           {+q\nu_1 \nu_2 \nu_3 y(z/q)=0,}$\\
         $q$-P$(E_3^{(1)})(b)$
           &$\displaystyle
           {z\left(z-\frac{\kappa_1}{\nu_8}\right)y(qz)
            +\Bigl( -{\nu_1}z^2+qdz-\frac{q\prod_{i=1}^5\nu_i}{\kappa_2}\Bigr)y(z)}$\\
           &$\displaystyle
           {+q\nu_1 \nu_2 \nu_3 (z-\nu_4)y(z/q)=0,}$\\
         $q$-P$(E_2^{(1)})(a)$
           &$\displaystyle
           {\frac{\kappa_1}{\nu_8} z (qz-\nu_4) y(qz)
    +\Bigl(\nu_1z^2+dz+\frac{q\prod_{i=1}^5\nu_i}{\kappa_2}\Bigr)y(z)
    +q\nu_1 \nu_2 \nu_3 y(z/q)=0,}$\\
         $q$-P$(E_2^{(1)})(b)$
           &$\displaystyle
           {z\left(z-\frac{\kappa_1}{\nu_8}\right)y(qz)
             +\left(-{\nu_1}z^2+dz\right)y(z)
             -q\nu_1 \nu_2 \nu_3 (z-\nu_4)y(z/q)=0,}$\\
         $q$-P$(A_1^{(1)})$
           &$\displaystyle
           {z\left(z-\frac{\kappa_1}{\nu_8}\right)y(qz)
             +dzy(z)
             -q\nu_1 \nu_2 \nu_3 (z-\nu_4)y(z/q)=0,}$\\
         $q$-P$\underset{|\alpha^2|=8}{(A_1^{(1)})}$
           &$\displaystyle
           {{\frac{\kappa_1}{\nu_8}}z(qz-\nu_4)y(qz)
    + (\nu_1z^2-dz ) y(z)
    +q\nu_1\nu_2\nu_3y(z/q)=0.}$
      \end{tabular}
  \endgroup}
\end{center}
The Newton diagrams of these equations are described as follows.
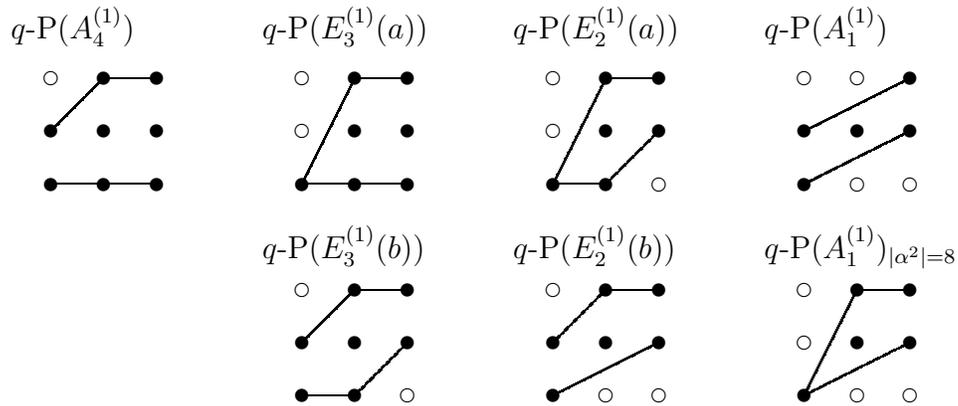
\begin{center}
\begin{picture}(370,175)(0,0)
  \put(15,155){$q$-P($A_4^{(1)}$)}
    \put(30,140){\circle{5}}
    \put(50,140){\circle*{5}}
    \put(70,140){\circle*{5}}
    \put(30,120){\circle*{5}}
    \put(50,120){\circle*{5}}
    \put(70,120){\circle*{5}}
    \put(30,100){\circle*{5}}
    \put(50,100){\circle*{5}}
    \put(70,100){\circle*{5}}
    \qbezier(50,140)(60,140)(70,140)
%    \qbezier(30,120)(30,130)(30,100)
    \qbezier(30,100)(50,100)(70,100)
    \qbezier(50,140)(40,130)(30,120)
%    \qbezier(70,100)(70,120)(70,140)
  \put(110,155){$q$-P($E_3^{(1)}(a)$)}
    \put(125,140){\circle{5}}
    \put(145,140){\circle*{5}}
    \put(165,140){\circle*{5}}
    \put(125,120){\circle{5}}
    \put(145,120){\circle*{5}}
    \put(165,120){\circle*{5}}
    \put(125,100){\circle*{5}}
    \put(145,100){\circle*{5}}
    \put(165,100){\circle*{5}}
    \qbezier(145,140)(155,140)(165,140)
%    \qbezier(165,140)(165,120)(165,100)
    \qbezier(165,100)(145,100)(125,100)
    \qbezier(125,100)(135,120)(145,140)
  \put(110,72){$q$-P($E_3^{(1)}(b)$)}
    %\put(90,110){\vector(1,-4){15}}
    \put(125,60){\circle{5}}
    \put(145,60){\circle*{5}}
    \put(165,60){\circle*{5}}
    \put(125,40){\circle*{5}}
    \put(145,40){\circle*{5}}
    \put(165,40){\circle*{5}}
    \put(125,20){\circle*{5}}
    \put(145,20){\circle*{5}}
    \put(165,20){\circle{5}}
    \qbezier(145,60)(155,60)(165,60)
%    \qbezier(165,60)(165,50)(165,40)
    \qbezier(165,40)(155,30)(145,20)
    \qbezier(125,40)(135,50)(145,60)
%    \qbezier(125,40)(125,30)(125,20)
    \qbezier(125,20)(135,20)(145,20)
  \put(205,155){$q$-P($E_2^{(1)}(a)$)}
    %\put(185,117.5){\vector(1,0){15}}
    \put(220,140){\circle{5}}
    \put(240,140){\circle*{5}}
    \put(260,140){\circle*{5}}
    \put(220,120){\circle{5}}
    \put(240,120){\circle*{5}}
    \put(260,120){\circle*{5}}
    \put(220,100){\circle*{5}}
    \put(240,100){\circle*{5}}
    \put(260,100){\circle{5}}
    \qbezier(240,140)(250,140)(260,140)
%    \qbezier(260,140)(260,130)(260,120)
    \qbezier(260,120)(250,110)(240,100)
    \qbezier(240,100)(230,100)(220,100)
    \qbezier(220,100)(230,120)(240,140)
  \put(205,72){$q$-P($E_2^{(1)}(b)$)}
    %\put(185,110){\vector(1,-4){15}}
    %\put(185,37.5){\vector(1,0){15}}
    \put(220,60){\circle{5}}
    \put(240,60){\circle*{5}}
    \put(260,60){\circle*{5}}
    \put(220,40){\circle*{5}}
    \put(240,40){\circle*{5}}
    \put(260,40){\circle*{5}}
    \put(220,20){\circle*{5}}
    \put(240,20){\circle{5}}
    \put(260,20){\circle{5}}
    \qbezier(240,60)(250,60)(260,60)
%    \qbezier(260,60)(260,50)(260,40)
    \qbezier(260,40)(240,30)(220,20)
%    \qbezier(220,20)(220,30)(220,40)
    \qbezier(220,40)(230,50)(240,60)
  \put(300,155){$q$-P($A_1^{(1)}$)}
    %\put(280,117.5){\vector(1,0){15}}
    \put(315,140){\circle{5}}
    \put(335,140){\circle{5}}
    \put(355,140){\circle*{5}}
    \put(315,120){\circle*{5}}
    \put(335,120){\circle*{5}}
    \put(355,120){\circle*{5}}
    \put(315,100){\circle*{5}}
    \put(335,100){\circle{5}}
    \put(355,100){\circle{5}}
%    \qbezier(355,140)(355,130)(355,120)
    \qbezier(355,140)(335,130)(315,120)
%    \qbezier(315,100)(315,110)(315,120)
    \qbezier(315,100)(335,110)(355,120)
  \put(300,72){$q$-P${(A_1^{(1)})_{|\alpha^2|=8}}$}
    %\put(280,110){\vector(1,-4){15}}
    %\put(280,37.5){\vector(1,0){15}}
    \put(315,60){\circle{5}}
    \put(335,60){\circle*{5}}
    \put(355,60){\circle*{5}}
    \put(315,40){\circle{5}}
    \put(335,40){\circle*{5}}
    \put(355,40){\circle*{5}}
    \put(315,20){\circle*{5}}
    \put(335,20){\circle{5}}
    \put(355,20){\circle{5}}
    \qbezier(335,60)(345,60)(355,60)
%    \qbezier(355,60)(355,50)(355,40)
    \qbezier(355,40)(335,30)(315,20)
    \qbezier(315,20)(325,40)(335,60)
  \end{picture}
\end{center}
In section \ref{sec:cqHs}, we introduce the confluent $q$-Heun equation, the biconfluent $q$-Heun equation and the doubly confluent $q$-Heun equation.
The equation corresponding to $q$-P$(A_4^{(1)})$ is the confluent $q$-Heun equation, the equation corresponding to $q$-P$(E_3^{(1)})(a)$ is the biconfluent $q$-Heun equation, and the equation corresponding to $q$-P$(E_3^{(1)})(b)$ is the doubly confluent $q$-Heun equation.
Moreover, the equation corresponding to $q$-P$(E_2^{(1)})(b)$ is the simply-reduced doubly confluent $q$-Heun equation, and the equation corresponding to $q$-P$(A_1^{(1)})$ is the doubly-reduced doubly confluent $q$-Heun equation.

\section{Confluent $q$-Heun equations} \label{sec:cqHs}

We define the confluent $q$-Heun equation by
\begin{align}
& (\alpha_1x+\alpha_0)f(qx) -(\beta_2x^2+\beta_1x+\beta_0)f(x) +(\gamma_2x^2+\gamma_1x+\gamma_0)f(x/q)=0 
\label{eq:cqHE}
\end{align}
under the condition $\alpha_1 \alpha_0\gamma_2\gamma_0\not=0$.
The equations
\begin{align}
& (\alpha_2x^2+\alpha_1x+\alpha_0)f(qx) -(\beta_2x^2+\beta_1x+\beta_0)f(x) +(\gamma_1x+\gamma_0)f(x/q)=0, \label{eq:cqHE2} \\
& (\alpha_2x^2+\alpha_1x+\alpha_0)f(qx) -(\beta_2x^2+\beta_1x+\beta_0)f(x) + (\gamma_2 x^2 +\gamma_1 x)f(x/q)=0, \label{eq:cqHE3} \\
& (\alpha_2x^2+\alpha_1x )f(qx) -(\beta_2x^2+\beta_1x+\beta_0)f(x) +(\gamma_2x^2+\gamma_1x+\gamma_0) f(x/q)=0 \label{eq:cqHE4}
\end{align}
are transformed to equation (\ref{eq:cqHE}) on other parameters by a gauge transformation and/or by replacing the variable as $x \to 1/x $.
We also call equation (\ref{eq:cqHE2}) (resp.~(\ref{eq:cqHE3}), (\ref{eq:cqHE4})) the confluent $q$-Heun equation under the condition $\alpha_2\alpha_0 \gamma_1 \gamma_0\not=0$ (resp.~$\alpha_2\alpha_0 \gamma_2 \gamma_1 \not=0$, $\alpha_2\alpha_1 \gamma_2 \gamma_0 \not=0$).
In equation (\ref{eq:cqHE}), if $\beta _2 =0$ (resp.~$\beta _2 \not = 0$), then we call it the reduced (resp.~non-reduced) confluent $q$-Heun equation.
Similarly, we can define the reduced or non-reduced confluent $q$-Heun equation for equations (\ref{eq:cqHE2}), (\ref{eq:cqHE3}) and (\ref{eq:cqHE4}).
The Newton polygons of the non-reduced confluent $q$-Heun equation in equations (\ref{eq:cqHE}), (\ref{eq:cqHE2}), (\ref{eq:cqHE3}) and (\ref{eq:cqHE4}) are described as
\begin{center}
\begin{picture}(370,85)(0,0)
    \put(20,60){\circle*{5}}
    \put(17,36){$*$}
%    \put(20,40){\circle*{5}}
    \put(20,20){\circle*{5}}
    \put(40,60){\circle*{5}}
    \put(37,36){$*$}
%    \put(40,40){\circle*{5}}    
    \put(37,16){$*$}
    \put(60,60){\circle{5}}
    \put(60,40){\circle*{5}}
%    \put(40,20){\circle{5}}
    \put(60,20){\circle*{5}}
    \qbezier(20,60)(30,60)(40,60)
    \qbezier(40,60)(50,50)(60,40)
    \qbezier(20,20)(40,20)(60,20)
    \put(160,60){\circle*{5}}
    \put(157,36){$*$}
%    \put(160,40){\circle*{5}}
    \put(160,20){\circle*{5}}
    \put(140,60){\circle*{5}}
    \put(137,36){$*$}
%    \put(140,40){\circle*{5}}    
    \put(137,16){$*$}
    \put(120,60){\circle{5}}
    \put(120,40){\circle*{5}}
%    \put(120,20){\circle{5}}
    \put(120,20){\circle*{5}}
    \qbezier(160,60)(150,60)(140,60)
    \qbezier(140,60)(130,50)(120,40)
    \qbezier(120,20)(140,20)(160,20)
    \put(260,60){\circle*{5}}
    \put(257,36){$*$}
%    \put(260,40){\circle*{5}}
    \put(260,20){\circle*{5}}
    \put(240,20){\circle*{5}}
    \put(237,36){$*$}
%    \put(240,40){\circle*{5}}    
    \put(237,56){$*$}
    \put(220,60){\circle*{5}}
    \put(220,40){\circle*{5}}
    \put(220,20){\circle{5}}
    \qbezier(260,60)(240,60)(220,60)
    \qbezier(240,20)(230,30)(220,40)
    \qbezier(240,20)(250,20)(260,20)
    \put(320,60){\circle*{5}}
    \put(317,36){$*$}
%    \put(320,40){\circle*{5}}
    \put(320,20){\circle*{5}}
    \put(340,20){\circle*{5}}
    \put(337,36){$*$}
%    \put(340,40){\circle*{5}}    
    \put(337,56){$*$}
    \put(360,60){\circle*{5}}
    \put(360,40){\circle*{5}}
    \put(360,20){\circle{5}}
    \qbezier(320,60)(340,60)(360,60)
    \qbezier(340,20)(350,30)(360,40)
    \qbezier(340,20)(330,20)(320,20)
  \end{picture}
\end{center}
The Newton polygons of the reduced confluent $q$-Heun equation in equations (\ref{eq:cqHE}), (\ref{eq:cqHE2}), (\ref{eq:cqHE3}) and (\ref{eq:cqHE4}) are described as
\begin{center}
\begin{picture}(370,85)(0,0)
    \put(20,60){\circle*{5}}
    \put(17,36){$*$}
%    \put(20,40){\circle*{5}}
    \put(20,20){\circle*{5}}
    \put(40,60){\circle{5}}
    \put(37,36){$*$}
%    \put(40,40){\circle*{5}}    
    \put(37,16){$*$}
    \put(60,60){\circle{5}}
    \put(60,40){\circle*{5}}
%    \put(40,20){\circle{5}}
    \put(60,20){\circle*{5}}
%    \qbezier(20,60)(30,60)(40,60)
    \qbezier(20,60)(40,50)(60,40)
    \qbezier(20,20)(40,20)(60,20)
    \put(160,60){\circle*{5}}
    \put(157,36){$*$}
%    \put(160,40){\circle*{5}}
    \put(160,20){\circle*{5}}
    \put(140,60){\circle{5}}
    \put(137,36){$*$}
%    \put(140,40){\circle*{5}}    
    \put(137,16){$*$}
    \put(120,60){\circle{5}}
    \put(120,40){\circle*{5}}
%    \put(120,20){\circle{5}}
    \put(120,20){\circle*{5}}
%    \qbezier(160,60)(150,60)(140,60)
    \qbezier(160,60)(140,50)(120,40)
    \qbezier(120,20)(140,20)(160,20)
    \put(260,60){\circle*{5}}
    \put(257,36){$*$}
%    \put(260,40){\circle*{5}}
    \put(260,20){\circle*{5}}
    \put(240,20){\circle{5}}
    \put(237,36){$*$}
%    \put(240,40){\circle*{5}}    
    \put(237,56){$*$}
    \put(220,60){\circle*{5}}
    \put(220,40){\circle*{5}}
    \put(220,20){\circle{5}}
    \qbezier(260,60)(240,60)(220,60)
    \qbezier(260,20)(240,30)(220,40)
%    \qbezier(240,20)(250,20)(260,20)
    \put(320,60){\circle*{5}}
    \put(317,36){$*$}
%    \put(320,40){\circle*{5}}
    \put(320,20){\circle*{5}}
    \put(340,20){\circle{5}}
    \put(337,36){$*$}
%    \put(340,40){\circle*{5}}    
    \put(337,56){$*$}
    \put(360,60){\circle*{5}}
    \put(360,40){\circle*{5}}
    \put(360,20){\circle{5}}
    \qbezier(320,60)(340,60)(360,60)
    \qbezier(320,20)(340,30)(360,40)
%    \qbezier(340,20)(330,20)(320,20)
  \end{picture}
\end{center}
We define the biconfluent $q$-Heun equation by
\begin{align}
& f(qx) -(\beta_2x^2+\beta_1x+\beta_0)f(x) +(\gamma_2x^2+\gamma_1x+\gamma_0)f(x/q)=0
\label{eq:bqHE}
\end{align}
under the condition $\beta_2 \gamma_2\gamma_0\not=0$.
The equations
\begin{align}
& (\alpha_2x^2+\alpha_1x+\alpha_0)f(qx) -(\beta_2x^2+\beta_1x+\beta_0)f(x) + f(x/q)=0, \label{eq:bqHE2} \\
& (\alpha_2x^2+\alpha_1x+\alpha_0)f(qx) -(\beta_2x^2+\beta_1x+\beta_0)f(x) + x^2 f(x/q)=0, \label{eq:bqHE3} \\
& x^2 f(qx) -(\beta_2x^2+\beta_1x+\beta_0)f(x) +(\gamma_2x^2+\gamma_1x+\gamma_0) f(x/q)=0 , \label{eq:bqHE4} \\
& (\alpha_1x+\alpha_0)f(qx) -(\beta_2x^2+\beta_1x+\beta_0)f(x) + (\gamma_1x+\gamma_0) f(x/q)=0, \label{eq:bqHE5}\\
& (\alpha_2x^2+\alpha_1x )f(qx) -(\beta_2x^2+\beta_1x+\beta_0)f(x) + (\gamma_2x^2+\gamma_1x ) f(x/q)=0 \label{eq:bqHE6} 
\end{align}
are transformed to equation (\ref{eq:bqHE}) on other parameters by a gauge transformation and/or by replacing the variable as $x \to 1/x $.
We also call equation (\ref{eq:bqHE2}) (resp.~(\ref{eq:bqHE3}), (\ref{eq:bqHE4}), (\ref{eq:bqHE5}), (\ref{eq:bqHE6})) the biconfluent $q$-Heun equation under the condition $\alpha_2\alpha_0 \beta_2 \not=0$ (resp.~$\alpha_2\alpha_0 \beta_0 \not=0$, $\beta_0 \gamma_2 \gamma_0 \not=0$, $\alpha_1\alpha_0 \beta_2 \gamma_1 \gamma_0 \not=0$, $\alpha_2\alpha_1 \beta_0 \gamma_2 \gamma_1 \not=0$).
The Newton polygons of the biconfluent $q$-Heun equation in equations (\ref{eq:bqHE}), (\ref{eq:bqHE2}), (\ref{eq:bqHE3}), (\ref{eq:bqHE4}), (\ref{eq:bqHE5}) and (\ref{eq:bqHE6}) are described as
\begin{center}
\begin{picture}(510,85)(0,0)
    \put(20,60){\circle*{5}}
    \put(17,36){$*$}
%    \put(20,40){\circle*{5}}
    \put(20,20){\circle*{5}}
    \put(40,60){\circle*{5}}
    \put(37,36){$*$}
%    \put(40,40){\circle*{5}}    
    \put(37,16){$*$}
    \put(60,60){\circle{5}}
    \put(60,40){\circle{5}}
%    \put(40,20){\circle{5}}
    \put(60,20){\circle*{5}}
    \qbezier(20,60)(30,60)(40,60)
    \qbezier(40,60)(50,40)(60,20)
    \qbezier(20,20)(40,20)(60,20)
 
   \put(130,60){\circle*{5}}
    \put(127,36){$*$}
%    \put(130,40){\circle*{5}}
    \put(130,20){\circle*{5}}
    \put(110,60){\circle*{5}}
    \put(107,36){$*$}
%    \put(110,40){\circle*{5}}    
    \put(107,16){$*$}
    \put(90,60){\circle{5}}
    \put(90,40){\circle{5}}
%    \put(90,20){\circle{5}}
    \put(90,20){\circle*{5}}
    \qbezier(130,60)(120,60)(110,60)
    \qbezier(110,60)(100,40)(90,20)
    \qbezier(90,20)(110,20)(130,20)  

   \put(200,60){\circle*{5}}
    \put(197,36){$*$}
%    \put(200,40){\circle*{5}}
    \put(200,20){\circle*{5}}
    \put(180,20){\circle*{5}}
    \put(177,36){$*$}
%    \put(190,40){\circle*{5}}    
    \put(177,56){$*$}
    \put(160,60){\circle*{5}}
    \put(160,40){\circle{5}}
%    \put(160,20){\circle{5}}
    \put(160,20){\circle{5}}
    \qbezier(200,60)(180,60)(160,60)
    \qbezier(160,60)(170,40)(180,20)
    \qbezier(180,20)(190,20)(200,20)  

   \put(230,60){\circle*{5}}
    \put(227,36){$*$}
%    \put(230,40){\circle*{5}}
    \put(230,20){\circle*{5}}
    \put(250,20){\circle*{5}}
    \put(247,36){$*$}
%    \put(260,40){\circle*{5}}    
    \put(247,56){$*$}
    \put(270,60){\circle*{5}}
    \put(270,40){\circle{5}}
%    \put(270,20){\circle{5}}
    \put(270,20){\circle{5}}
    \qbezier(270,60)(250,60)(230,60)
    \qbezier(270,60)(260,40)(250,20)
    \qbezier(250,20)(240,20)(230,20)  

    \put(300,60){\circle{5}}
    \put(320,60){\circle*{5}}
    \put(340,60){\circle{5}}
    \put(300,40){\circle*{5}}
    \put(317,36){$*$}
%    \put(320,40){\circle*{5}}
    \put(340,40){\circle*{5}}
    \put(300,20){\circle*{5}}
    \put(317,16){$*$}
%    \put(320,20){\circle{5}}
    \put(340,20){\circle*{5}}
    \qbezier(300,40)(310,50)(320,60)
    \qbezier(320,60)(330,50)(340,40)
    \qbezier(300,20)(320,20)(340,20)

    \put(370,20){\circle{5}}
    \put(390,20){\circle*{5}}
    \put(410,20){\circle{5}}
    \put(370,40){\circle*{5}}
    \put(387,36){$*$}
%    \put(390,40){\circle*{5}}
    \put(410,40){\circle*{5}}
    \put(370,60){\circle*{5}}
    \put(387,56){$*$}
%    \put(390,20){\circle{5}}
    \put(410,60){\circle*{5}}
    \qbezier(370,40)(380,30)(390,20)
    \qbezier(390,20)(400,30)(410,40)
    \qbezier(370,60)(390,60)(410,60)
%    \qbezier(410,40)(410,30)(410,20)
%    \qbezier(410,20)(390,30)(370,40)
%    \qbezier(370,40)(370,50)(370,60)
  \end{picture}
\end{center}
We define the doubly confluent $q$-Heun equation by
\begin{align}
& x f(qx) -(\beta_2x^2+\beta_1x+\beta_0)f(x) +(\gamma_2x^2+\gamma_1x+\gamma_0)f(x/q)=0, 
\label{eq:dqHE}
\end{align}
under the condition $\gamma_2\gamma_0\not=0$.
The equations
\begin{align}
& (\alpha_2x^2+\alpha_1x+\alpha_0)f(qx) -(\beta_2x^2+\beta_1x+\beta_0)f(x) + x f(x/q)=0, \label{eq:dqHE2} \\
& (\alpha_1x+\alpha_0)f(qx) -(\beta_2x^2+\beta_1x+\beta_0)f(x) + (\gamma_2x ^2+\gamma_1x ) f(x/q)=0, \label{eq:dqHE3} \\
& (\alpha_2 x^2 +\alpha_1 x)f(qx) -(\beta_2x^2+\beta_1x+\beta_0)f(x) + (\gamma_1 x +\gamma_0 ) f(x/q)=0 \label{eq:dqHE4}
\end{align}
are transformed to equation (\ref{eq:dqHE}) on other parameters by a gauge transformation and/or by replacing the variable as $x \to 1/x $.
We also call equation (\ref{eq:dqHE2}) (resp.~(\ref{eq:dqHE3}), (\ref{eq:dqHE4})) the doubly confluent $q$-Heun equation under the condition $\alpha_2\alpha_0 \not=0$ (resp.~$\alpha_1\alpha_0 \gamma_2 \gamma_1 \not=0$, $\alpha_2\alpha_1 \gamma_1 \gamma_0 \not=0$).
In equation (\ref{eq:dqHE}), if $\beta _2 =0$ or $\beta _0 =0$  (resp.~$\beta _2 \not = 0 \not  = 0 \beta$), then we call it the reduced (resp.~non-reduced) doubly confluent $q$-Heun equation.
If $\beta _2 = \beta _0 =0$, then we call it the doubly-reduced doubly confluent $q$-Heun equation.
Similar definitions hold for equations (\ref{eq:dqHE2}), (\ref{eq:dqHE3}) and (\ref{eq:dqHE4}).
The Newton polygons of the non-reduced doubly confluent $q$-Heun equation in equations (\ref{eq:dqHE}), (\ref{eq:dqHE2}), (\ref{eq:dqHE3}) and (\ref{eq:dqHE4}) are described as
\begin{center}
\begin{picture}(370,85)(0,0)
    \put(20,60){\circle*{5}}
    \put(17,36){$*$}
%    \put(20,40){\circle*{5}}
    \put(20,20){\circle*{5}}
    \put(40,60){\circle*{5}}
    \put(37,36){$*$}
%    \put(40,40){\circle*{5}}    
%    \put(37,16){$*$}
    \put(40,20){\circle*{5}}
    \put(60,60){\circle{5}}
    \put(60,40){\circle*{5}}
%    \put(40,20){\circle{5}}
    \put(60,20){\circle{5}}
    \qbezier(20,60)(30,60)(40,60)
    \qbezier(40,60)(50,50)(60,40)
    \qbezier(40,20)(50,30)(60,40)
    \qbezier(20,20)(30,20)(40,20)
    \put(160,60){\circle*{5}}
    \put(157,36){$*$}
%    \put(160,40){\circle*{5}}
    \put(160,20){\circle*{5}}
    \put(140,60){\circle*{5}}
    \put(137,36){$*$}
    \put(140,20){\circle*{5}}    
%    \put(137,16){$*$}
    \put(120,60){\circle{5}}
    \put(120,40){\circle*{5}}
%    \put(120,20){\circle{5}}
    \put(120,20){\circle{5}}
    \qbezier(160,60)(150,60)(140,60)
    \qbezier(140,60)(130,50)(120,40)
    \qbezier(140,20)(130,30)(120,40)
    \qbezier(140,20)(150,20)(160,20)
    \put(260,60){\circle{5}}
%    \put(257,36){$*$}
    \put(260,40){\circle*{5}}
    \put(260,20){\circle*{5}}
    \put(240,20){\circle*{5}}
    \put(237,36){$*$}
    \put(240,60){\circle*{5}}    
%    \put(237,56){$*$}
    \put(220,60){\circle*{5}}
    \put(220,40){\circle*{5}}
    \put(220,20){\circle{5}}
    \qbezier(220,60)(230,60)(240,60)
    \qbezier(260,40)(250,50)(240,60)
    \qbezier(240,20)(230,30)(220,40)
    \qbezier(240,20)(250,20)(260,20)
    \put(320,60){\circle{5}}
%    \put(317,36){$*$}
    \put(320,40){\circle*{5}}
    \put(320,20){\circle*{5}}
    \put(340,20){\circle*{5}}
    \put(337,36){$*$}
    \put(340,60){\circle*{5}}    
%    \put(337,56){$*$}
    \put(360,60){\circle*{5}}
    \put(360,40){\circle*{5}}
    \put(360,20){\circle{5}}
    \qbezier(340,60)(350,60)(360,60)
    \qbezier(340,60)(330,50)(320,40)
    \qbezier(340,20)(350,30)(360,40)
    \qbezier(340,20)(330,20)(320,20)
  \end{picture}
\end{center}
The Newton polygons of the two singly-reduced confluent $q$-Heun equation and the doubly-reduced confluent $q$-Heun equation in equation (\ref{eq:dqHE}) are described as
\begin{center}
%\begin{picture}(0,85)(0,0)
\begin{picture}(250,85)(0,0)
    \put(20,60){\circle*{5}}
    \put(17,36){$*$}
%    \put(20,40){\circle*{5}}
    \put(20,20){\circle*{5}}
    \put(40,60){\circle{5}}
    \put(37,36){$*$}
%    \put(40,40){\circle*{5}}    
%    \put(37,16){$*$}
    \put(40,20){\circle*{5}}
    \put(60,60){\circle{5}}
    \put(60,40){\circle*{5}}
%    \put(40,20){\circle{5}}
    \put(60,20){\circle{5}}
%    \qbezier(20,60)(30,60)(40,60)
    \qbezier(20,60)(40,50)(60,40)
    \qbezier(40,20)(50,30)(60,40)
    \qbezier(20,20)(30,20)(40,20)

    \put(120,60){\circle*{5}}
    \put(117,36){$*$}
%    \put(120,40){\circle*{5}}
    \put(120,20){\circle*{5}}
    \put(140,60){\circle*{5}}
    \put(137,36){$*$}
%    \put(140,40){\circle*{5}}    
%    \put(137,16){$*$}
    \put(140,20){\circle{5}}
    \put(160,60){\circle{5}}
    \put(160,40){\circle*{5}}
%    \put(140,20){\circle{5}}
    \put(160,20){\circle{5}}
    \qbezier(120,60)(130,60)(140,60)
    \qbezier(140,60)(150,50)(160,40)
    \qbezier(120,20)(140,30)(160,40)
%    \qbezier(120,20)(130,20)(140,20)

    \put(220,60){\circle*{5}}
    \put(217,36){$*$}
%    \put(220,40){\circle*{5}}
    \put(220,20){\circle*{5}}
    \put(240,60){\circle{5}}
    \put(237,36){$*$}
%    \put(240,40){\circle*{5}}    
%    \put(237,16){$*$}
    \put(240,20){\circle{5}}
    \put(260,60){\circle{5}}
    \put(260,40){\circle*{5}}
%    \put(240,20){\circle{5}}
    \put(260,20){\circle{5}}
%    \qbezier(20,60)(30,60)(40,60)
    \qbezier(220,60)(240,50)(260,40)
    \qbezier(220,20)(240,30)(260,40)
%    \qbezier(220,20)(230,20)(240,20)
  \end{picture}
\end{center}
The Newton polygons of the two singly-reduced confluent $q$-Heun equation and the doubly-reduced confluent $q$-Heun equation in equation (\ref{eq:dqHE4}) are described as
\begin{center}
%\begin{picture}(0,85)(0,0)
\begin{picture}(250,85)(0,0)
    \put(20,60){\circle{5}}
%    \put(17,36){$*$}
    \put(20,40){\circle*{5}}
    \put(20,20){\circle*{5}}
    \put(40,20){\circle*{5}}
    \put(37,36){$*$}
    \put(40,60){\circle{5}}    
%    \put(37,56){$*$}
    \put(60,60){\circle*{5}}
    \put(60,40){\circle*{5}}
    \put(60,20){\circle{5}}
    \qbezier(60,60)(40,50)(20,40)
    \qbezier(40,20)(50,30)(60,40)
    \qbezier(40,20)(30,20)(20,20)

    \put(120,60){\circle{5}}
%    \put(117,36){$*$}
    \put(120,40){\circle*{5}}
    \put(120,20){\circle*{5}}
    \put(140,20){\circle{5}}
    \put(137,36){$*$}
    \put(140,60){\circle*{5}}    
%    \put(137,56){$*$}
    \put(160,60){\circle*{5}}
    \put(160,40){\circle*{5}}
    \put(160,20){\circle{5}}
    \qbezier(140,60)(150,60)(160,60)
    \qbezier(140,60)(130,50)(120,40)
    \qbezier(120,20)(140,30)(160,40)

    \put(220,60){\circle{5}}
%    \put(217,36){$*$}
    \put(220,40){\circle*{5}}
    \put(220,20){\circle*{5}}
    \put(240,20){\circle{5}}
    \put(237,36){$*$}
    \put(240,60){\circle{5}}    
%    \put(237,56){$*$}
    \put(260,60){\circle*{5}}
    \put(260,40){\circle*{5}}
    \put(260,20){\circle{5}}
    \qbezier(260,60)(240,50)(220,40)
    \qbezier(220,20)(240,30)(260,40)

  \end{picture}
\end{center}
The Newton polygons of the singly-reduced confluent $q$-Heun equations and the doubly-reduced confluent $q$-Heun equations in equations (\ref{eq:dqHE2}) and (\ref{eq:dqHE3}) are described similarly.

We investigate the limit procedure from the $q$-Heun equation or its confluent equations to the differential equations.
Let us consider the $q$-difference equation
\begin{align}
& (a_{-,2} x^2 + a_{-,1} x +a_{-,0} ) g(x/q) + (a_{+,2} x^2 + a_{+,1} x +a_{+,0} ) g(q x)  \label{eq:qHEa-a+} \\
& \qquad  + (a_{0,2} x^2 + a_{0,1} x +a_{0,0} ) g(x) = 0 . \nonumber
\end{align}
Set $q=1+ \varepsilon $.
By using Taylor's expansion
\begin{align}
& g(x/q) =g(x) + (-\varepsilon +\varepsilon ^2 )xg'(x) + \varepsilon ^2 x^2 g''(x) /2 +O( \varepsilon ^3),\\
& g(qx) =g(x) + \varepsilon  x g'(x) + \varepsilon ^2 x^2 g''(x) /2 +O( \varepsilon ^3), \nonumber 
\end{align}
we have
%\begin{align}
%& (a_{-,2} x^2 + a_{-,1} x +a_{-,0} ) (g(x) + (-\varepsilon +\varepsilon ^2 )xg'(x) + \varepsilon ^2 x^2 g''(x) /2 +O( \varepsilon ^3))  \\
%& +  (a_{+,2} x^2 + a_{+,1} x +a_{+,0} ) (g(x) + \varepsilon  x g'(x) + \varepsilon ^2 x^2 g''(x) /2 +O( \varepsilon ^3)) \nonumber \\
%& \qquad  + (a_{0,2} x^2 + a_{0,1} x +a_{0,0} ) g(x) = 0 . \nonumber
%\end{align}
\begin{align}
& \{ (a_{-,2} + a_{+,2} +a_{0,2} ) x^2 + (a_{-,1} + a_{+,1} +a_{0,1} ) x +(a_{-,0} + a_{+,0} +a_{0,0} ) \}  g(x) \\
& +  \varepsilon  \{ (a_{+,2} - a_{-,2}) x^2 + ( a_{+,1} - a_{-,1}) x + (a_{+,0}- a_{-,0}) \} x g'(x) \nonumber \\
& + \varepsilon ^2  \{ (a_{+,2} + a_{-,2}) x^2 + (a_{+,1} + a_{-,1}) x + (a_{+,0} + a_{-,0}) \} x^2 g''(x) /2  \nonumber \\
& + (a_{-,2} x^2 + a_{-,1} x +a_{-,0} )  \varepsilon ^2  xg'(x) +O( \varepsilon ^3) =0 . \nonumber 
\end{align}
If there exist constants $b_0, b_1, b_2, b_{0,1} , b_{1,1} , b_{2,1} , b_{0,0} , b_{1,0} , b_{2,0} $ such that
%\tilde{b}_0 , \tilde{b}_1 , \tilde{b}_2 , \tilde{\tilde{b}}_0 , \tilde{\tilde{b}}_1 , \tilde{\tilde{b}}_2 $ such that
\begin{align}
& (a_{-,2} + a_{+,2} +a_{0,2} ) / \varepsilon ^2 \to  b_{2,0} , \;  (a_{-,1} + a_{+,1} +a_{0,1} ) / \varepsilon ^2 \to  b_{1,0} , \\
&  (a_{-,0} + a_{+,0} +a_{0,0} ) / \varepsilon ^2 \to  b_{0,0} , \nonumber \\
& (a_{+,2} - a_{-,2}) / \varepsilon \to  b_{2,1} , \;  (a_{+,1} - a_{-,1}) / \varepsilon \to  b_{1,1} , \;  (a_{+,0} - a_{-,0}) / \varepsilon \to  b_{0,1} , \nonumber \\ 
& (a_{+,2} + a_{-,2}) /2 \to  b_{2} , \;  (a_{+,1} + a_{-,1}) / 2 \to  b_{1} , \;  (a_{+,0} + a_{-,0}) / 2 \to  b_{0} \nonumber 
\end{align}
as $\varepsilon \to 0$, then 
\begin{align}
& a_{+,2} \to b_2 ,\; a_{-,2} \to  b_{2} , \; a_{0,2} \to -2 b_2, \; a_{+,1} \to b_1 ,\; a_{-,1} \to b_{1} , \; a_{0,1} \to -2 b_1, \\
& a_{+,0} \to b_0 ,\; a_{-,0} \to  b_{0} , \; a_{0,0} \to -2 b_0 \nonumber 
\end{align}
as $\varepsilon \to 0$, and we obtain the differential equation
\begin{align}
&  x ^2 ( b_2 x^2 + b_1 x + b_0 ) g''(x) + x \{ (b_{2,1} +b _2 ) x^2 + (b_{1,1} +b_1) x + (b_{0,1}+b_0 ) \} g'(x) \label{eq:DEep0} \\
&  \qquad + (b_{2,0} x^2 + b_{1,0} x + b_{0,0}  ) g(x) =0 . \nonumber 
\end{align}

If $b_{0,0}=0$, $b_2\neq 0 \neq b_0$ and $b_1 ^2 \neq 4 b_0 b_2$ in equation (\ref{eq:DEep0}), then the differential equation is written as
\begin{align}
&   b_2 x ( x- \alpha _1 ) ( x- \alpha _2 ) g''(x) + \{ (b_{2,1} +b _2 ) x^2 + (b_{1,1} +b_1) x + (b_{0,1}+b_0 ) \} g'(x) \\
&  \qquad + (b_{2,0} x + b_{1,0}  ) g(x) =0 ,  \nonumber 
\end{align}
where $\{ \alpha _1 ,\alpha _2 \}$ is the set of the roots of the quadratic equation $ b_2 x^2 + b_1 x + b_0 =0 $, and we obtain the Heun equation (HE) in equation (\ref{eq:Heun}) by changing the variable as $x \to \alpha _1 x$.
If $b_{0,0} \neq 0$, $b_2\neq 0 \neq b_0$ and $b_1 ^2 \neq 4 b_0 b_2$ in equation (\ref{eq:DEep0}), then we obtain the Heun equation by using the gauge-transformation $g(x) \to x^{\rho } g(x) $ for some $\rho $ and changing the variable $x$.

If $b_{0,0}=0$, $b_2= 0 $ and $b_1 \neq 0 \neq  b_0$ in equation (\ref{eq:DEep0}), then the differential equation is written as
\begin{align}
&  x  ( b_1 x + b_0 ) g''(x) + \{ b_{2,1} x^2 + (b_{1,1} +b_1) x + (b_{0,1}+b_0 ) \} g'(x) \\
&  \qquad + (b_{2,0} x + b_{1,0}  ) g(x) =0 ,  \nonumber 
\end{align}
and we obtain the confluent Heun equation (CHE) in equation (\ref{eq:CHE}) by changing the variable $x $.
If $b_{0,0} \neq 0$, $b_2= 0 $ and $b_1 \neq 0 \neq  b_0$ in equation (\ref{eq:DEep0}), then we obtain the confluent Heun equation by using the gauge-transformation $g(x) \to x^{\rho } g(x) $ for some $\rho $ and changing the variable $x $.
If $a_{+,2} =0 $ or $a_{-,2} =0 $ in equation (\ref{eq:qHEa-a+}) (i.e.~the case of the confluent $q$-Heun equation (\ref{eq:cqHE})), then we have $b_2 =0$, and the corresponding differential equation is CHE or more degenerated ones.
%On the confluent $q$-Heun equation (\ref{eq:cqHE}), the condition $a_{+,2} =0 $ is required.
If $a_{0,2} =0 $ in addition, then we have $(a_{-,2} + a_{+,2} ) / \varepsilon  \to 0$ as $\varepsilon \to 0$, and it follows from $a_{+,2} =0 $ or $a_{-,2} =0 $ that $(a_{-,2} - a_{+,2} ) / \varepsilon  \to 0$ as $\varepsilon \to 0$.
Then we have $b_2 = b_{2,1} =0$, and the corresponding differential equation is the reduced CHE or more degenerated ones, where the irregular singularity $z= \infty $ is ramified.

If $b_{0,0}=0$, $b_2= b_1 =0 $ and $b_0 \neq 0 $ in equation (\ref{eq:DEep0}), then the differential equation is written as
\begin{align}
&  b_0 x g''(x) + \{ b_{2,1} x^2 + b_{1,1} x + (b_{0,1}+b_0 ) \} g'(x) + (b_{2,0} x + b_{1,0}  ) g(x) =0 ,
\end{align}
and we obtain the biconfluent Heun equation (BHE) in equation (\ref{eq:BHE}).
If $b_{0,0} \neq 0$, $b_2= b_1 =0 $ and $b_0 \neq 0 $ in equation (\ref{eq:DEep0}), then we obtain the biconfluent Heun equation by using the gauge-transformation $g(x) \to x^{\rho } g(x) $ for some $\rho $.
If $a_{+,2} = a_{+,1} = 0 $ in equation (\ref{eq:qHEa-a+}) (i.e.~the case of the biconfluent $q$-Heun equation (\ref{eq:bqHE})), then we have $b_2 =b_1= 0$, and the corresponding differential equation is BHE or more degenerated ones.

If $b_{0,0}=0$, $b_2= b_0 =0 $ and $b_1 \neq 0 $ in equation (\ref{eq:DEep0}), then the differential equation is written as
\begin{align}
&  b_1 x^2 g''(x) + \{ b_{2,1} x^2 + (b_{1,1} +b_1) x + b_{0,1} \} g'(x)  + (b_{2,0} x + b_{1,0}  ) g(x) =0 ,
\end{align}
and we obtain the doubly confluent Heun equation (DHE) in equation (\ref{eq:DHE}).
If $b_{0,0} \neq 0$, $b_2= b_0 =0 $ and $b_1 \neq 0 $ and in equation (\ref{eq:DEep0}), then we obtain the doubly confluent Heun equation by using the gauge-transformation $g(x) \to x^{\rho } g(x) $ for some $\rho $.
% $b_{0,1} \neq 0 $
If $a_{+,2} = a_{+,0} = 0 $ in equation (\ref{eq:qHEa-a+}) (i.e.~the case of the doubly confluent $q$-Heun equation (\ref{eq:dqHE})), then we have $b_2= b_0 =0 $, and the corresponding differential equation is DHE or more degenerated ones.
If $a_{0,2} =0 $ (resp.~$a_{0,0} =0 $) in addition, then we have $b_2 = b_{2,1} =0$ (resp.~$b_0 = b_{0,1} =0$), which is obtained as the case of the reduced CHE.
The corresponding differential equation is the reduced DHE or more degenerated ones, where the irregular singularity $z= \infty $ (resp.~the origin $z=0$) is ramified.
If $a_{+,2} = a_{+,0} = a_{0,2} =a_{0,0} =0 $, then we have $b_2 = b_{2,1} =b_0 = b_{0,1} =0$, and the corresponding differential equation is doubly-reduced DHE or more degenerated ones.
%If $b_{0,0} \neq 0$, $b_2= b_0 =0 $, $b_1 \neq 0 $ and $b_{0,1} = 0 $ in equation (\ref{eq:DEep0}), then the differential equation is written as
%\begin{align}
%&  b_1 x^3 g''(x) +  x^2 \{ (b_{2,1} +b _2 ) x + (b_{1,1} +b_1)  \} g'(x) + (b_{2,0} x^2 + b_{1,0} x + b_{0,0}  ) g(x) =0 ,
%\end{align}
%which corresponds to the reduced doubly confluent Heun equation (RDHE) which has a ramified irregular singularity about $x=0$.

\section{Concluding remarks}  \label{sec:rmk}

In this paper, we obtained several degenerated $q$-Heun equations from Lax forms associated to several $q$-Painlev\'e equations, and we used the Lax forms of Murata and the ones of Kajiwara, Noumu and Yamada.
By picking up some degenerated $q$-Heun equations obtained from the $q$-Painlev\'e equations, we defined the confluent $q$-Heun equation, the biconfluent $q$-Heun equation, the doubly confluent $q$-Heun equation and the reduced ones, and investigated the limit to the differential equations.
Consequently, we obtaind the confluent Heun equation, the biconfluent Heun equation, the doubly confluent Heun equation and the reduced ones by the limit $q\to 1$.

We present another method to guess the limit $q \to 1$ by considering local solutions at $x=0$ and $x=\infty $.
Recall the $q$-difference equation 
\begin{align}
&  (qx-a_1t)f(qx)-\{q^2\kappa_1x^2+dx+(\theta_1+\theta_2)t\}f(x) \label{eq:qDEMA4} \\
&  +\kappa_1\kappa_2(qx-a_3)(x-a_2t)f(x/q)=0, \nonumber
\end{align}
which was obtained from a piece of Lax form of type $A_4$ of Murata.
We substitute the expansion
\begin{equation}
f(x)=x^\rho(c_0+c_1x+c_2x^2+\cdots), \; c_0\not=0
\label{eq:expa0}
\end{equation}
at $x=0$ to the $q$-difference equation.
Then, we obtain
\begin{equation}
  a_1q^{2\rho}+(\theta_1+\theta_2)q^\rho-\kappa_1\kappa_2a_2a_3=0
\label{eq:qrhoM}
\end{equation}
by observing the coefficient of $x^{\rho }$.
Hence, we obtain two values $q^{\rho }$ from equation (\ref{eq:qrhoM}), and the expansions in equation (\ref{eq:expa0}) for two values $\rho$ are jist like the ones for the regular singularity of the second order linear differential equation.
On the other hand, we  substitute the expansion
\begin{equation}
f(x)=x^{-\rho}(c_0+c_1x^{-1}+c_2x^{-2}+\cdots), \; c_0\not=0
\end{equation}
at $x=\infty $ to the $q$-difference equation.
Then, we obtain
\begin{equation}
 \kappa_2q^{\rho}-q=0
\end{equation}
by observing the coefficient of $x^{2-\rho}$.
We obtain one value $q^{\rho} = q / \kappa_2$, and the situation is different from the regular singularity of the second order linear differential equation.
Therefore, we detect that the corresponging differential equation has a regular singularity at $x=0$ and an irregular singularity at $x=\infty $.
Note that this observation is related to the shape of the Newton diagram associated to the $q$-difference equation (\ref{eq:qDEMA4}) (see section \ref{sec:Mlist}).
We might investigate the differential equations which are related to the $q$-difference equation obtained from the types $A_6$ and $A_7$ of Murata and some types of Kajiwara, Noumi and Yamada.

\section*{Acknowledgements}
The second author was supported by JSPS KAKENHI Grant Number JP22K03368.


\begin{thebibliography}{9999}
%\bibitem{Adams28}
%C.~R.~Adams, On the linear ordinary q-difference equation. {\it Ann. Math.} {\bf 30} (1928), 195--205.
\bibitem{Adams31}
C.~R.~Adams, Linear $q$-difference equations. {\it Bull. Amer. Math. Soc.} {\bf 37} (1931), 361--400.
%Adams, C. R. (1931). Linear $ q $-difference equations. Bull. Amer. Math. Soc., 37(12), 361-400.
\bibitem{Hahn}
W.~Hahn, On linear geometric difference equations with accessory parameters, {\it Funkcial. Ekvac.} {\bf 14} (1971), 73--78.
\bibitem{Htd}
Y.~Hatsuda, Quasinormal modes of Kerr-de Sitter black holes via the Heun function, \textit{Class. Quantum Grav.} \textbf{38} (2021), 025015.
%Yasuyuki Hatsuda
%DOI 10.1088/1361-6382/abc82e
%\bibitem{HMST}
%N.~Hatano, R.~Matsunawa, T.~Sato, K.~Takemura, Variants of $q$-hypergeometric equation, to appear in {\it Funkcial. Ekvac.}, {\sf arXiv:1910.12560}.
\bibitem{JS} M.~Jimbo, H.~Sakai, A $q$-analog of the sixth Painlev\'e equation. {\it Lett. Math. Phys.}  {\bf 38} (1996), 145--154
\bibitem{KNY} K.~Kajiwara, M.~Noumi, Y.~Yamada, Geometric aspects of Painlev\'e equations. {\it J. Phys. A} {\bf 50} (2017), 073001.
\bibitem{MST}
R.~Matsunawa, T.~Sato, K.~Takemura, Variants of confluent $q$-hypergeometric equations, in Geometric and Harmonic Analysis on Homogeneous Spaces and Applications, Springer Proceedings in Mathematics \& Statistics {\bf 366}, 161--180, Springer, Cham, 2021. {\sf arXiv:2005.13223}.
%Matsunawa R., Sato T., Takemura K. (2021) Variants of Confluent q-Hypergeometric Equations. In: Baklouti A., Ishi H. (eds) Geometric and Harmonic Analysis on Homogeneous Spaces and Applications. TJC 2019. Springer Proceedings in Mathematics & Statistics, vol 366. Springer, Cham. pp 161-180, https://doi.org/10.1007/978-3-030-78346-4_10
\bibitem{Mrt}
M.~Murata, Lax forms of the $q$-Painlev\'e equations. {\it J. Phys. A} {\bf 42} (2009), 115201, 17 pp.
\bibitem{NM}
S.~Noda, H.~Motohashi, Spectroscopy of $\mathrm{Kerr}\text{\ensuremath{-}}{\mathrm{AdS}}_{5}$ spacetime with the Heun function: Quasinormal modes, greybody factor, and evaporation, \textit{Phys. Rev. D} \textbf{106} (2022), 064025, 22 pp.
%  author = {Noda, Sousuke and Motohashi, Hayato},
%  doi = {10.1103/PhysRevD.106.064025},
\bibitem{Ohy}
Y.~Ohyama, A unified approach to $q$-special functions of the Laplace type, {\sf arXiv:1103.5232}.
%\bibitem{Ron}
%A. Ronveaux A. (ed.), Heun's differential equations, Oxford Science Publications, Oxford University Press, Oxford, 1995.
\bibitem{RuiN} S.~N.~M.~Ruijsenaars, Integrable $BC_N$ analytic difference operators: Hidden parameter symmetries and eigenfunctions, in Proceedings Cadiz 2002 NATO Advanced Research Workshop "New trends in integrability and partial solvability", NATO Science Series {\bf 132}, 217--261, Kluwer, Dordrecht, 2004.
%\bibitem{RuiN} Ruijsenaars S.~N.~M.~, Integrable $BC_N$ analytic difference operators: Hidden parameter symmetries and eigenfunctions, NATO Science Series {\bf 132}, 217--261, Kluwer, Dordrecht, 2004.
%\bibitem{Rui06} S.~N.~M.~Ruijsenaars, Eigenfunctions with a zero eigenvalue for differences of elliptic relativistic Calogero-Moser Hamiltonians, {\it Theoret. and Math. Phys.} \textbf{146} (2006), 25--33.
%\bibitem{Rui09} S.~N.~M.~Ruijsenaars, Hilbert-Schmidt operators vs. integrable systems of elliptic Calogero-Moser type. I. The eigenfunction identities, {\it Comm. Math. Phys.} \textbf{286} (2009), 629--657.
\bibitem{Sak} 
H.~Sakai, Rational surfaces associated with affine root systems and geometry of the Painlev\'e equations, {\it Commun. Math. Phys.}, {\bf 220} (2001), 165--229.
%\bibitem{SY} 
%H.~Sakai, M.~Yamaguchi, Spectral types of linear $q$-difference equations and $q$-analog of middle convolution, {\textit{Int. Math. Res. Not.}} \textbf{2017} (2017), 1975--2013.
\bibitem{STT1}
S.~Sasaki, S.~Takagi, K.~Takemura, $q$-Heun equation and initial-value space of $q$-Painlev\'e equation, {\it Contemp. Math.} {\bf 782} (2023), 119--142.
\bibitem{STT2}
S.~Sasaki, S.~Takagi, K.~Takemura, $q$-middle convolution and $q$-Painlev\'e equation, {\textit{SIGMA}} \textbf{18} (2022), paper 056.
%, {\sf arXiv:2201.03960}.
\bibitem{Sau}
J.~Sauloy, Basic Modern Theory of Linear Complex Analytic $q$-Difference Equations. Vol. 287. American Mathematical Society, 2024.
%Sauloy, Jacques. Basic Modern Theory of Linear Complex Analytic $q$-Difference Equations. Vol. 287. American Mathematical Society, 2024.
%\bibitem{SL}
%S. Slavyanov and W. Lay, Special Functions., Oxford Science Publications, Oxford University Press, Oxford, 2000.
%\bibitem{Tak2}
%K. Takemura, The Heun equation and the Calogero-Moser-Sutherland system {\rm II}: the perturbation and the algebraic solution, Electron. J. Differential Equations, {\bf 2004} (2004), no. 15 1--30. 
\bibitem{TakHP} K.~Takemura, Heun equation and Painlev\'e equation, in Proceedings of Workshop on Elliptic Integrable Systems 2004 Kyoto, Rokko Lectures in Mathematics {\bf 18}, 305--322, 2005.
% (math.CA/0503288) 
\bibitem{TakM}
K.~Takemura, Middle convolution and Heun's equation, {\it SIGMA} {\bf 5} (2009), paper 040.
%\bibitem{TakS} K.~Takemura, Heun's differential equation, Selected papers on analysis and differential equations, 45--68, Amer. Math. Soc. Transl. Ser. 2, {\bf 230}, Amer. Math. Soc., Providence, 2010.
%\bibitem{TakIT} 
%K.~Takemura, Integral transformation of Heun's equation and some applications, {\it J. Math. Soc. Japan} {\bf 69} (2017), 849--891.
\bibitem{TakR}
K.~Takemura, Degenerations of Ruijsenaars-van Diejen operator and $q$-Painleve equations, {\it J. Integrable Systems} {\bf 2} (2017), xyx008.
%, {\sf arXiv:1608.07265}.
%\bibitem{TakqH}
%K.~Takemura., On $q$-deformations of the Heun equation, {\it SIGMA} {\bf 14} (2018), paper 061.
\bibitem{vD0} J.~F.~van Diejen, Integrability of difference Calogero-Moser systems. {\it J. Math. Phys.} {\bf 35} (1994), 2983--3004. 
%\bibitem{WW}
%E.~T.~Whittaker, G.~N.~Watson, A course of modern analysis. Fourth edition. Cambridge University Press, New York, 1962.
\bibitem{Ye} Y.~Yamada, A Lax formalism for the elliptic difference Painlev\'e equation, {\it SIGMA}, {\bf 5} (2009), paper 042.
\bibitem{Y} Y.~Yamada, Lax formalism for $q$-Painleve equations with affine Weyl group symmetry of type $E^{(1)}_n$, {\it Int. Math. Res. Notices}, {\bf 2011} (2011), 3823--3838.


\end{thebibliography}
\end{document}